\documentclass[12pt,leqno]{article}

\usepackage{amsthm}
\usepackage[usenames]{color}

\usepackage{amsmath}
\usepackage{amscd}
\usepackage{amsmath}
\usepackage{amssymb}
\usepackage{latexsym}
\makeatletter

\@addtoreset{equation}{section}
\makeatother

\newtheorem{thm}{Theorem}[subsection]
\newtheorem{pr}[thm]{Proposition}
\newtheorem{df}[thm]{Definition}
\newtheorem{lm}[thm]{Lemma}
\newtheorem{cor}[thm]{Corollary}

\newtheorem{rmk}[thm]{Remark}
\newtheorem{ex}[thm]{Example}

\newcommand{\sm}{\raisebox{2.33pt}{~\rule{6.4pt}{1.3pt}~}}

\input xy \xyoption{all} \CompileMatrices

\begin{document}

\title{Ramification groups of
coverings and valuations} 

\author{Takeshi Saito}
\maketitle

\begin{abstract}
We give a purely scheme theoretic
construction of the filtration
by ramification groups of
the Galois group of a covering.
The valuation need not be discrete
but the normalizations are required
to be locally of complete intersection.
\end{abstract}

For a Galois extension of a complete
discrete valuation field with
not necessarily perfect residue field,
the filtration by ramification groups
on the Galois group is defined in
a joint article \cite{AS} with Ahmed Abbes.
Although the definition in \cite{AS} is based on
rigid geometry, it is later observed
that the use of rigid geometry can be
avoided and the conventional
language of schemes suffices
(\cite{mu}).
In this article, we reformulate the construction in
\cite{AS} in the language of schemes.
As a byproduct, we give a generalization
for ramified finite Galois coverings 
of normal noetherian schemes
and valuations not necessarily discrete.

All the ideas are present in
\cite{AS}, possibly in different formulation.
As in \cite{AS},
the main ingredients in the definition
of ramification groups are the followings:
First, we interpret a subgroup as a
quotient of the fiber functor
with a cocartesian property,
Proposition \ref{prGN}.
Thus, the definition of ramification groups
is a consequence of a construction
of quotients of the fiber functor,
indexed by elements of the rational
value group of valuation.

The required quotients of the fiber functor
are constructed as the sets of connected
components of geometric fibers of
dilatations defined by
an immersion of the covering to
a smooth scheme over the base scheme.
Here a crucial ingredient is
the reduced fiber theorem of
Bosch-L\"utkebohmert-Raynaud \cite{BLR}
recalled in Theorem \ref{thmRFT}.
This specializes to
the finiteness theorem of
Grauert-Remmert 
in the classical case where
the base is a discrete valuation ring.
A variant of the filtration 
is defined using the underlying sets of
geometric fibers of quasi-finite
schemes without using
the sets of connected
components.

To prove basic properties of
ramification groups stated in
Theorem \ref{thmgr}
including the rationality of breaks,
semi-continuity etc.,
a key ingredient is a generalization
due to Temkin \cite{Tst} of 
the semi-stable reduction theorem of
curves recalled in Theorem \ref{thmsst}.

Let $X$ be a normal noetherian scheme
and $U\subset X$ be a dense open 
subscheme.
The Zariski-Riemann space 
$\tilde X$ is defined
as the inverse limit of
proper schemes $X'$ over $X$
such that $U'=U\times_XX'
\to U$ is an isomorphism.
Points of $\tilde X$ on the boundary
$\tilde X\sm U$ correspond
bijectively to the inverse limits of the images 
of the closed points
by the liftings of the morphisms $T={\rm Spec}\, A
\to X$ for valuation rings
$A\subsetneqq K=k(t)$ for 
points $t\in U$ such that
$T\times_XU$ consists of the single points $t$.

Let $W\to U$ be a finite \'etale connected
Galois covering of Galois group $G$.
We will construct in Theorem \ref{thmgr}
filtrations 
$(G_T^\gamma)$ and 
$(G_T^{\gamma+})$ on $G$
by ramification groups for a morphism
$T\to X$ as above
indexed by the positive part 
$(0,\infty)_{\Gamma_{\mathbf Q}}
\subset \Gamma_{\mathbf Q}
=\Gamma\otimes {\mathbf Q}$
for the value group
$\Gamma=K^\times/A^\times$.
To complete the definition,
we need to assume that
for every intermediate covering $V\to U$,
the normalization $Y$ of $X$
in $V$ is locally of complete intersection
over $X$ to assure
the cocartesian property in
Proposition \ref{prGN}.
The required cocartesian property
Proposition \ref{prcoca}
is then a consequence of
a lifting property in commutative
algebra recalled in
Proposition \ref{prCA}.

The definition depends on $X$
not only on $W\to U$.
In other words,
for a normal noetherian scheme
$X'$ over $X$ as above,
the filtrations 
$(G_T^\gamma)$ and 
$(G_T^{\gamma+})$ defined for $X$
and those for $X'$ may be different.
This arises from the fact that
the formation of the normalization $Y$
need not commute with base change
$X'\to X$.
To obtain a definition depending only on
$W\to U$, one would need
to take inverse limit
with respect to $X'$. This requires
that the normalizations over
$T$ to be locally of complete intersection.

By Proposition \ref{prGN},
the definition of the filtrations
$(G_T^\gamma)$ and 
$(G_T^{\gamma+})$ are reduced to
the construction of surjections
$F_T^\infty \to F_T^\gamma$
and
$F_T^\infty \to F_T^{\gamma+}$
for a fiber functor $F_T^\infty$.
To define them,
for each intermediate covering
$V\to U$, we take an embedding
$Y\to Q$ of the normalization
to a smooth scheme over $X$.
Further taking a ramified covering and
a blow-up $X'$, we find an effective
Cartier divisor $R'\subset X'$
and a lifting $T'\to X'$
of $T\to X$ such that the valuation
$v'(R')$ of $R'$ is $\gamma$
for each
$\gamma\in \Gamma_{\mathbf Q}$.
Then, we define a dilatation
$Q^{\prime (R')}$ over $X'$
to be the normalization of an open
subscheme $Q^{\prime [R']}$
of the blow-up of the base change
$Q'=Q\times_XX'$ at the closed
subscheme
$Y\times_XR'\subset Q\times_XX'$.
To obtain a construction independent of
the choice of $X'$,
we apply the reduced fiber theorem
of Bosch-L\"utkebohmert-Raynaud 
for $Q^{\prime (R')}\to X'$
to be flat and to have
reduced geometric fibers.

Now the desired functor
$F_T^\gamma(Y/X)$
is defined as the set of connected
components of the geometric fiber
of $Q^{\prime (R')}\to X'$
at the image of the closed point by
$T'\to X'$.
Example \ref{exQ}.1 and Remark \ref{rmk}
imply that we recover the construction
in \cite{AS} in the classical case
where $X=T$ is the spectrum of
a complete discrete valuation ring.
Its variant $F_T^{\gamma+}(Y/X)$
is defined more simply as
the geometric fiber of the inverse
image $Y'\times_{Q^{\prime [R']}}
Q^{\prime (R')}$ with respect to the
morphism
$Y'=Y\times_XX'\to Q^{\prime [R']}$
lifting the original immersion $Y\to Q$.
The fact that the construction is
independent of the choice of
immersion $Y\to Q$ is based
on a homotopy invariance
of dilatations proved in Proposition
\ref{prFunQ}.

To study the behavior of the functors
$F_T^\gamma$ and
$F_T^{\gamma+}$ thus defined
for variable $\gamma$,
we use a semi-stable curve
$C$ over $X$ defined by
$st=f$ for a non-zero divisor 
$f$ on $X$ defining
an effective Cartier divisor
$D\subset X$ such that
$D\cap U=\varnothing$ as
a parameter space for $\gamma$.
Let $\tilde D\subset C$ denote
the effective Cartier divisor defined by $t$.
Then, for $\gamma\in [0,v(D)]_{\Gamma_
{\mathbf Q}}$,
there is a lifting $T'\to C$ of
$T\to X$ such that
$v'(\tilde D)=\gamma$.
Using this together with 
a local description 
Proposition \ref{prDiv} of Cartier divisors
on a semi-stable curve over
a normal noetherian scheme
and 
a combination of the reduced fiber theorem
and the semi-stable reduction theorem
over a general base scheme,
we derive basic properties
of $F_T^\gamma$ and
$F_T^{\gamma+}$ in Proposition \ref{prstD}
to prove Theorem \ref{thmval}
and Theorem \ref{thmgr}.

The author heartily thanks Ahmed Abbes
for collaboration on ramification theory
and encouragement on writing this article.
This work is supported by JSPS Grants-in-Aid 
for Scientific Research
(A) 26247002.

\subsection*{Convention}
In this article, we assume that
for a noetherian scheme $X$,
the normalization of a
scheme of finite type over $X$
remains to be of finite type over $X$.
This property is satisfied
if $X$ is of finite type
over a field, ${\mathbf Z}$
or a complete discrete valuation ring,
for example.

\tableofcontents
\section{Preliminaries}

\subsection{Connected components}


\begin{df}[{\rm \cite[D\'efinition (6.8.1)]{EGA4}}]\label{dfsep}
Let $f\colon X\to S$ be a flat morphism
locally of finite presentation 
of schemes.
We say that $f$ is {\em reduced}
if for every geometric point $s$ of $S$,
the geometric fiber $X_s$ is {\em reduced}.
\end{df}

In SGA 1 Expos\'e X Definition 1.1,
reduced morphism is called separable morphism.

We study the sets of connected components
of geometric fibers of a flat and reduced morphism
of finite type.
Let $S$ be a scheme
and $s$ and $t$ be geometric points
of $S$.
Let $S_{(s)}$
denote the strict localization.
A specialization $s\gets t$
of geometric points
means a morphism
$S_{(s)}\gets t$ over $S$.

Assume that $S$ is noetherian.
Let $X\to S$ be a flat and reduced morphism
of finite type
and let $s\gets t$ be a specialization
of geometric points of $S$.
We define the cospecialization
mapping 
\begin{equation}
\pi_0(X_s)\to \pi_0(X_t)
\label{eqcosp}
\end{equation}
as follows.
By replacing $S$ by 
the closure of the image of $t$,
we may assume that $S$
is integral and
that $t$ is above the generic point
$\eta$ of $S$.
By replacing $S$ further
by a quasi-finite scheme over $S$
such that the function field is
a finite extension of
$\kappa(\eta)$ in $\kappa(t)$,
we may assume that
the canonical mapping
$\pi_0(X_t)\to \pi_0(X_\eta)$
is a bijection.
Let $U\subset S$ be a dense open
subset such that the canonical mapping
$\pi_0(X_\eta)\to \pi_0(X_U)$
is a bijection.
Then, by \cite[Corollaire (18.9.11)]{EGA4},
the canonical mapping
$\pi_0(X_U)\to \pi_0(X)$
is also a bijection.
Thus, we define the cospecialization
mapping (\ref{eqcosp})
to be the composition
$$\pi_0(X_s)\to \pi_0(X)
\overset\simeq\gets
\pi_0(X_\eta)
\overset\simeq\gets \pi_0(X_t).$$

We say that the sets of connected components of geometric
fibers of $X\to S$
are locally constant
if for every specialization
$s\gets t$ of geometric points of $S$,
the cospecialization mapping
$\pi_0(X_s)\to \pi_0(X_t)$
is a bijection.
By \cite[Th\'eor\`eme (9.7.7)]{EGA4}
and by noetherian induction,
there exists a finite
stratification $S=\coprod_iS_i$
by locally closed subschemes
such that the sets of connected components of geometric
fibers of the base change
$X\times_SS_i\to S_i$
are locally constant
for every $i$.
We call this fact that
the sets of connected components of geometric
fibers of $X\to S$
are constructible.

\begin{rmk}\label{rmk}
{\rm Let $S={\rm Spec}\, {\cal O}_K$
for a discrete valuation ring
${\cal O}_K$ and $X={\rm Spec}\, A$
be an affine scheme of finite type over $S$.
Let $\bar s\to S$ be the geometric closed point.
Let ${\mathfrak X}={\rm Spf}\, \hat A$
be the formal completion along the closed fiber
and ${\mathfrak X}_{\bar K}={\rm Sp}\, 
\hat A\otimes_{{\cal O}_K}\! \bar K$
be the associated affinoid variety
over an algebraic closure $\bar K$
of the fraction field $K$ of ${\cal O}_K$.
If $X$ is flat and reduced over $S$,
the cospecialization mapping
$\pi_0(X_{\bar s})
\to \pi_0({\mathfrak X}_{\bar K})$
is a bijection.}
\end{rmk}

Let $Y\to S$ be another flat and reduced morphism
of finite type
and let $f\colon X\to Y$
be a morphism over $S$.
The cospecialization
mappings (\ref{eqcosp})
form a commutative diagram
\begin{equation}
\begin{CD}
\pi_0(X_s)@>>> \pi_0(X_t)\\
@VVV@VVV\\
\pi_0(Y_s)@>>> \pi_0(Y_t).
\end{CD}
\label{eqcosp2}
\end{equation}

\begin{lm}\label{lmA}
Let $f\colon X\to Y$ be a morphism
of schemes of finite type
over a noetherian scheme $S$.
Assume that $X$ is \'etale over $S$
and that $Y$ is flat and reduced over $S$.
Let $A$ denote the subset
of $X$ consisting of the
images of geometric points
$x$ of $X$ satisfying the following
conditions:
Let $s$ be the geometric point of 
$S$ defined as the image of $x$
and let $C\subset Y_s$
be the connected component
of the fiber
containing the image of $x$.
Then, $f_s^{-1}(C)\subset X_s$
consists of a single point $x$.

Then $A$ is closed.
\end{lm}

\proof{
By the constructibility of
connected components of
geometric fibers of $Y$,
the subset $A\subset X$
is constructible.
For a specialization
$s\gets t$ of geometric points
of $S$, the upper horizontal arrow
in the commutative diagram
$$\begin{CD}
X_s@>>>X_t\\
@VVV@VVV\\
 \pi_0(Y_s)@>>> \pi_0(Y_t)
\end{CD}$$
is an injection since $X\to S$ is \'etale.
Hence $A$ is closed under specialization and is closed.
\qed }

\medskip

We have specialization mappings
going the other way for 
proper morphisms.
Let $X$ be a proper scheme 
over $S$.
Let $s\gets t$ be a specialization
of geometric points of $S$.
Then, the inclusion
$X_s\to X\times_SS_{(s)}$
induces a bijection
$\pi_0(X_s)\to \pi_0(X\times_SS_{(s)})$
by \cite[IV Proposition (2.1)]{Arcata}.
Its composition with the mapping
$\pi_0(X_t)\to \pi_0(X\times_SS_{(s)})$
induced by the morphism
$X_t\to X\times_SS_{(s)}$
defines the specialization mapping
\begin{equation}
\pi_0(X_s)\gets \pi_0(X_t).
\label{eqsp}
\end{equation}
For a morphism $X\to Y$
of proper schemes over $S$,
the specialization mappings
make a commutative diagram
\begin{equation}
\begin{CD}
\pi_0(X_s)@<<< \pi_0(X_t)\\
@VVV@VVV\\
\pi_0(Y_s)@<<< \pi_0(Y_t).
\end{CD}
\label{eqspXY}
\end{equation}

\begin{lm}\label{lmB}
Let $f\colon X\to Y$ be a finite unramified morphism
of schemes.
Let $B$ denote the subset
of $X$ consisting of the
images of geometric points
$x$ of $X$ satisfying the following
conditions:
For the geometric point $y$ of 
$Y$ defined as the image of $x$,
the fiber $X\times_Yy$
consists of a single point $x$.

Then, $B$ is open.
\end{lm}

\proof{
The complement $X\sm B$
equals the image of the complement
$X\times_YX\sm X$
of the diagonal by a projection.
Since $X\to Y$ is unramified,
the complement
$X\times_YX\sm X
\subset
X\times_YX$ is closed.
Since the projection
$X\times_YX\to X$
is finite,
the image $X\sm B$ is closed.
\qed }

\begin{pr}\label{prCA}
Let 
\begin{equation}
\begin{CD}
Z'@>>>X'\\
@VVV
\hspace{-10mm}
\square
\hspace{7mm}
@VVfV
\\
Z@>>>X
\end{CD}
\label{eqCA}
\end{equation}
be a cartesian diagram
of noetherian schemes.
Assume that $X$ is normal,
the horizontal arrows
are closed immersion,
the right vertical arrow
is quasi-finite 
and the left vertical arrow is
finite.
Assume further 
that there exists a dense
open subscheme $U\subset X$
such that $U'=U\times_XX'\to U$ is faithfully flat
and that $U'\subset X'$ is also dense.

{\rm 1.}
Let $C\subset Z$
be an irreducible closed subset
and let $C'\subset f^{-1}(C)$ be an irreducible component.
Then, $C'\to C$ is surjective.

{\rm 2.}
Let $C\subset Z$
be a connected closed subset
and let $C'\subset f^{-1}(C)$ be a connected component.
Then, $C'\to C$ is surjective.
\end{pr}

\proof{
1.
By shrinking $U$ if necessary,
we may assume that $U'\to U$
is finite.
By Zariski's main theorem,
there exists a scheme
$\bar X'$ finite over $X$ containing
$X'$ as an open subscheme.
By replacing $\bar X'$
by the closure of $U'$,
we may assume that
$U'$ is dense in 
$\bar X'$.
Since $U'$ is closed in
$\bar X'\times_XU$,
we have 
$\bar X'\times_XU=U'$.
Since $Z'=(\bar X'\times_XZ)\cap X'$ is 
closed and open in
$\bar X'\times_XZ$,
by replacing $X'$ by $\bar X'$,
we may assume that $f$ is finite.

Since $f$ is a closed mapping,
it suffices to show
that the generic point
$z$ of $C$ is the image of 
the generic point $z'$ of $C'$.
Let $x'$ be a point of $C'$.
Replacing $X$ by an affine
neighborhood of $x=f(x')\in C$,
we may assume $X={\rm Spec}\ A$ and 
$X'={\rm Spec}\ B$
are affine.
Then, the assumption implies
that $A\to B$ is an injection
and $B$ is finite over $A$.
Since $x$ is a point of
the closure $C=\overline{\{z\}}$,
the assertion follows from
\cite[Chap.\ V, Section 2.4, Theorem 3]{Bo}.

2.
Let $C_1\subset C$ be an irreducible component
such that $C_1\cap f(C')$ is not empty.
Then, there exists an irreducible component
$C'_1$ of $f^{-1}(C_1)\subset f^{-1}(C)$ such that
$C'_1\cap C'$ is not empty.
By 1,
we have $C_1=f(C'_1)$.
Since $C'$ is a connected component
of $f^{-1}(C)$
and $C'_1\cap C'\neq \varnothing$,
we have $C'_1\subset C'$
and hence $C_1=f(C'_1)\subset f(C')$.
Thus, the complement $C\sm f(C')$
is the union of irreducible components
of $C$ not meeting $f(C')$
and is closed. 
Since $f(C')\subset C$ is also closed
and is non-empty,
we have $C=f(C')$.
\qed }

\begin{cor}\label{corCA}
Let 
\begin{equation}
\begin{CD}
Z'@>>>X'@<<<Y'_1@<{g'}<< Y'\\
@VVV
\hspace{-10mm}
\square
\hspace{7mm}
@VfVV
\hspace{-10mm}
\square
\hspace{7mm}
@V{f_1}VV@VV{f'}V
\\
Z@>>>X@<<<Y_1@<g<< Y
\end{CD}
\label{eqXYU}
\end{equation}
be a commutative diagram
of noetherian schemes
such that the left square is cartesian
and satisfies the conditions
in Proposition {\rm \ref{prCA}}.
Assume that $Y_1\subset X$ is a closed subscheme,
that the middle square is cartesian
and that the four arrows in
the right square are finite.
Assume that there exists a dense open
subscheme $V_1\subset Y_1$
such that 
$V=V_1\times_{Y_1}Y\subset Y$
is also dense and that 
$g|_V\colon V\to V_1$
and 
$g'|_{V'}\colon V'=V\times_YY'\to 
V'_1=V_1\times_{Y_1}Y'_1$
are isomorphisms.

{\rm 1.}
For any irreducible 
(resp.\ connected) component $C$
of $Y$,
we have
$f^{-1}(g(C))=g'(f^{\prime-1}(C))$.
Consequently, we have
$f^{-1}(g(Y))=g'(Y')$.

{\rm 2.}
Suppose that the mapping
$Z\times_XY\to\pi_0(Y)$ is a bijection.
Then, the diagram
\begin{equation}
\begin{CD}
Z'\cap Y'_1@<<< Z'\times_{X'} Y'\\
@VVV@VVV\\
Z\cap Y_1@<<< Z\times_XY
\end{CD}
\label{eqcoca1}
\end{equation}
of underlying sets
induces a surjection
$Z'\times_{X'} Y'
\to (Z'\cap Y'_1)
\times_{Z\cap Y_1}(Z\times_XY)$
of sets.
If $Z\times_XY\to Z\cap Y_1$
is surjective,
then $Z'\times_{X'} Y'\to Z'\cap Y'_1$
is also surjective.
Further if
$Y'\to Y$ is surjective,
then $Z'\cap Y'_1\to Z\cap Y_1$
is also surjective and
the diagram {\rm (\ref{eqcoca1})}
is a cocartesian diagram of underlying sets.

{\rm 3.}
The diagram
\begin{equation}
\begin{CD}
\pi_0(Z')@<<< Z'\cap Y'_1\\
@VVV@VVV\\
\pi_0(Z)@<<< Z\cap Y_1
\end{CD}
\label{eqcocaZ}
\end{equation}
of sets induces a surjection
$Z'\cap Y'_1\to \pi_0(Z')
\times_{\pi_0(Z)}(Z\cap Y_1)$ of sets.
If $Z\cap Y_1\to\pi_0(Z)$
is surjective, then
$Z'\cap Y'_1\to \pi_0(Z')$ is also surjective.
Further if $Z'\cap Y'_1\to Z\cap Y_1$ is surjective,
the diagram {\rm (\ref{eqcocaZ})}
is a cocartesian diagram of sets.
\end{cor}

\proof{
1.
Let $C\subset Y$ be an irreducible component.
The inclusion $f^{-1}(g(C))\supset g'(f^{\prime-1}(C))$
is clear. 
We show the other inclusion.
Since $V$ is dense in $Y$,
the intersection $C\cap V$
and hence its image
$g(C)\cap V_1$ are not empty.
Let $C'$ be an irreducible component of
$f^{-1}(g(C))\subset Y'_1$.
Since $Y'_1\to Y_1$ is finite
and $g(C)\subset Y_1$
is an irreducible closed subset,
we have $g(C)=f(C')$
by Proposition \ref{prCA}.1.
Since $f(C'\cap V'_1)
=f(C')\cap V_1
=g(C)\cap V_1$ is not empty,
$C'\cap V'_1
=g'(g^{\prime-1}(C'\cap V'_1))$
is also non-empty and 
hence is dense in $C'$.
Since 
$g^{\prime-1}(C'\cap V'_1)
= g^{\prime-1}(C')\cap V'
\subset f^{\prime-1}(C)$
and since $g'\colon Y'\to X'$ is proper,
we have $C'\subset g'(f^{\prime-1}(C))$.

Since a connected component of $Y$
and $Y$ itself are unions of irreducible
components of $Y$,
the remaining assertions follow
from the assertion for irreducible components.

2.
Let $z'\in Z'\cap Y'_1$ and
$y\in Z\times_XY$ be
points satisfying
$f(z')=g(y)$
in $Z\cap Y_1$.
Let $C\subset Y$ be the unique connected component
containing $y$.
Since $z'\in f^{-1}(g(C))
=g'(f^{\prime -1}(C))$ by 1,
there exists a point
$y'\in Z'\times_{X'}f^{\prime -1}(C)
\subset Z'\times_{X'}Y'$
such that
$z'=g'(y')$.
Since $f'(y')\in Z\times_XY$ 
is a unique point contained
in $C\in \pi_0(Y)$,
we have $y=f'(y')$.
Thus,
$(z',y)\in (Z'\cap Y'_1)
\times_{Z\cap Y_1}(Z\times_XY)$
is the image of 
$y'\in Z'\times_{X'}Y'$.

If $Z\times_XY\to Z\cap Y_1$
is surjective, then
$(Z'\cap Y'_1)
\times_{Z\cap Y_1}(Z\times_XY)\to Z'\cap Y'_1$ is surjective
and hence the first assertion implies
the surjectivity of
$Z'\times_{X'}Y'\to Z'\cap Y'_1$.

If both
$Z\times_XY\to Z\cap Y_1$
and $Y'\to Y$ are surjective,
then $Z'\times_{X'}Y'
=Z\times_XY'\to Z\times_XY$
is also surjective and hence
by the commutative diagram
(\ref{eqcoca1}),
the mapping 
$Z'\cap Y'_1\to Z\cap Y_1$
is a surjection.
This implies that the diagram
(\ref{eqcoca1}) with $Z'\times_{X'}Y'$
replaced by $(Z'\cap Y'_1)
\times_{Z\cap Y_1}(Z\times_XY)$
is a cocartesian diagram of 
underlying sets.
Hence the surjectivity of $Z'\times_{X'}Y'\to
Z'\cap Y'_1
\times_{Z\cap Y_1}(Z\times_XY)$
implies that the diagram (\ref{eqcoca1})
is a cocartesian diagram of 
underlying sets.

3.
Let $C'\subset Z'$ be a connected 
component 
and $z\in Z\cap Y_1$ be a point such that the connected component
$C\subset Z$ satisfying
$f(C')\subset C$ 
contains $z$.
Since $f(C')=C$ by Proposition \ref{prCA}.2,
the intersection
$C'\cap f^{-1}(z)
\subset Z'\cap Y'_1$
is not empty.
Hence $(C',z)\in
\pi_0(Z')
\times_{\pi_0(Z)}(Z\cap Y_1)$
is in the image of 
$C'\cap f^{-1}(z)
\subset Z'\cap Y'_1$.

The remaining assertions are proved
similarly as in 2. 
\qed }

\subsection{Flat and reduced morphisms}

Let $k\geqq 0$ be an integer.
Recall that a noetherian scheme
$X$ satisfies the condition {\rm (R$_k$)}
if for every point $x\in X$
of $\dim {\cal O}_{X,x}\leqq k$,
the local ring ${\cal O}_{X,x}$ is
regular
\cite[D\'efinition (5.8.2)]{EGA4}.
Recall also that a noetherian scheme
$X$ satisfies the condition {\rm (S$_k$)}
if for every point $x\in X$,
we have
${\rm prof}\, 
{\cal O}_{X,x}\geqq \inf(k,\dim {\cal O}_{X,x})$
\cite[D\'efinition (5.7.2)]{EGA4}.

\begin{pr}\label{prRk}
Let $f\colon X\to S$
be a flat morphism of finite
type of noetherian schemes
and $k\geqq 0$
be an integer.
We define a function $k\colon S\to{\mathbf N}$
by $k(s)=\max(k-\dim {\cal O}_{S,s},0)$.

{\rm 1.}
If $S$ satisfies the condition
{\rm (R$_k$)}
and if the fiber $X_s=X\times_Ss$
satisfies {\rm (R$_{k(s)}$)}
for every $s\in S$,
then $X$ satisfies the condition
{\rm (R$_k$)}.

{\rm 2.}
If $X$ satisfies the condition
{\rm (R$_k$)}
and if $f\colon X\to S$
is faithfully flat,
then $S$ satisfies the condition
{\rm (R$_k$)}.
\end{pr}

\proof{
1. Assume $\dim {\cal O}_{X,x}\leqq k$
and set $s=f(x)$.
Then, we have
$\dim {\cal O}_{S,s}\leqq 
\dim {\cal O}_{X,x}\leqq k$
and 
$\dim {\cal O}_{X_s,x}
=\dim {\cal O}_{X,x}
-\dim {\cal O}_{S,s}
\leqq k(s)$
by \cite[Proposition (5.7.2)]{EGA4}.
Hence ${\cal O}_{S,s}$
and ${\cal O}_{X_s,x}$
are regular by the assumption.
Thus 
${\cal O}_{X,x}$
is regular 
by \cite[Chap.\ 0$_{\rm IV}$ Proposition (17.3.3) (ii)]{EGA4}.

2. 
It follows from
\cite[Proposition (6.5.3) (i)]{EGA4}.
\qed }

\begin{pr}\label{prSk}
Let $f\colon X\to S$
be a flat morphism of finite
type of noetherian schemes
and $k\geqq 0$
be an integer.
Let the function $k\colon S\to{\mathbf N}$
be as in Proposition {\rm \ref{prRk}}.

{\rm 1.}
If $S$ satisfies the condition
{\rm (S$_k$)}
and if the fiber $X_s$
satisfies {\rm (S$_{k(s)}$)}
for every $s\in S$,
then $X$ satisfies the condition
{\rm (S$_k$)}.

{\rm 2.}
If $X$ satisfies the condition
{\rm (S$_k$)}
and if $f\colon X\to S$
is faithfully flat,
then $S$ satisfies the condition
{\rm (S$_k$)}.

{\rm 3.}
If $X$ satisfies the condition
{\rm (S$_k$)}
and if $S$ is of Cohen-Macaulay,
then the fiber $X_s$
satisfies {\rm (S$_{k(s)}$)}
for every $s\in S$.
\end{pr}

\proof{
1. Let $x\in X$
and $s=f(x)$.
Then, we have
${\rm prof}\, {\cal O}_{S,s}\geqq 
\inf (k, \dim {\cal O}_{S,s})$
and 
${\rm prof}\, {\cal O}_{X_s,x}\geqq 
\inf (k(s), \dim {\cal O}_{X_s,x})$
by the assumption.
By
$\dim {\cal O}_{X_s,x}
=\dim {\cal O}_{X,x}
-\dim {\cal O}_{S,s}$ 
\cite[Proposition (5.7.2)]{EGA4},
we have
$\inf (k, \dim {\cal O}_{S,s})
+
\inf (k(s), \dim {\cal O}_{X_s,x})
=\inf (k, \dim {\cal O}_{X,x})$.
Hence it follows from
${\rm prof}\, {\cal O}_{X,x}=
{\rm prof}\, {\cal O}_{S,s}
+{\rm prof}\, {\cal O}_{X_s,x}$
\cite[Proposition (6.3.1)]{EGA4}.

2. 
It follows from
\cite[Proposition (6.4.1) (i)]{EGA4}.

3. Let $x\in X$
and $s=f(x)$.
Then by the assumption, we have
${\rm prof}\, {\cal O}_{X,x}
\geqq
\inf (k, \dim {\cal O}_{X,x})$
and 
${\rm prof}\, {\cal O}_{S,s}
=\dim {\cal O}_{S,s}$.
By 
${\rm prof}\, {\cal O}_{X_s,x}=
{\rm prof}\, {\cal O}_{X,x}-
{\rm prof}\, {\cal O}_{S,s}
\geqq 0$
\cite[Proposition (6.3.1)]{EGA4}
and
$\dim {\cal O}_{X_s,x}=
\dim {\cal O}_{X,x}-
\dim {\cal O}_{S,s}$
\cite[Proposition (5.7.2)]{EGA4}
we have
${\rm prof}\, {\cal O}_{X_s,x}
\geqq
\inf (k-\dim {\cal O}_{S,s}, 
\dim {\cal O}_{X_s,x})\geqq k(s)$
and the assertion follows.
\qed }

\begin{cor}\label{corSk}
Let $f\colon X\to S$
be a flat morphism of finite
type of noetherian schemes
and let $U\subset X$
be the largest open subset
smooth over $S$.

{\rm 1.}
Assume that 
the fiber $X_s$ is reduced
for every $s\in S$.
Assume further that 
$S$ is normal and that 
for the generic point 
$s$ of each irreducible component,
$X_s$ is normal.
Then $X$ is normal.

{\rm 2.}
For $s\in S$ and a geometric point $\bar s$ 
above $s$, 
we consider the following conditions:

{\rm (1)} The geometric fiber
$X_{\bar s}$ is reduced.

{\rm (2)}
$U_s$ is dense in $X_s$.

Then, we have 
{\rm (1)}$\Rightarrow${\rm (2)}.
Conversely, if $X$ is normal and
$S$ is regular of dimension $\leqq 1$,
then we have
{\rm (2)}$\Rightarrow${\rm (1)}.
\end{cor}

\proof{
1. By Serre's criterion 
\cite[Th\'eor\`eme (5.8.6)]{EGA4},
$S$ satisfies {\rm (R$_2$)}
and {\rm (S$_1$)}.
By \cite[Proposition (5.8.5)]{EGA4},
every fiber $X_s$ satisfies {\rm (R$_1$)}
and {\rm (S$_0$)}.
Further if $s$ is 
the generic point of an irreducible component,
the fiber $X_s$ satisfies {\rm (R$_2$)}
and {\rm (S$_1$)}.
Since the function $k(s)$
for $k=2$ satisfies
$k(s)\leqq 1$ unless $s$ is 
the generic point an irreducible component
and $k(s)=2$ for such point,
the scheme $X$ satisfies the conditions  {\rm (R$_2$)}
and {\rm (S$_1$)}
by Propositions \ref{prRk}.1 and \ref{prSk}.1.
Thus the assertion follows
by \cite[Th\'eor\`eme (5.8.6)]{EGA4}.

2.
(1)$\Rightarrow$(2):
Since $X_{\bar s}$ is reduced,
there exists a dense open subset
$V\subset X_{\bar s}$
smooth over $\bar s$.
Since $f$ is flat,
the image of $V$ in $X_s$
is a subset of $U_s$.

(2)$\Rightarrow$(1):
Since $X$ satisfies
(S$_2$)
and $S$ is Cohen-Macaulay
of dimension $\leqq 1$,
the fiber $X_s$ satisfies
(S$_1$)
by Proposition \ref{prSk}.3.
Hence
the geometric fiber
$X_{\bar s}$
also satisfies
(S$_1$)
by \cite[Proposition (6.7.7)]{EGA4}.
By (2),
$X_{\bar s}$
satisfies
(R$_0$).
Hence the assertion follow from
\cite[Proposition (5.8.5)]{EGA4}.
\qed }

\begin{lm}\label{lmeta}
Let $S$ be a noetherian scheme
and let $f\colon Y\to X$ be a 
quasi-finite morphism of schemes of
finite type over $S$.
Assume that $X$ is smooth
over $S$ and that
$Y$ is flat and reduced over $S$.
Assume that there exist
dense open subschemes 
$U\subset S$
and $U\times_SX\subset
W\subset X$
such that
$Y\times_XW\to W$ is \'etale
and that for every point $s\in S$,
the inverse image 
$f_s^{-1}(W_s)
\subset Y_s=Y\times_Ss$
of $W_s=W\times_Ss
\subset X_s=X\times_Ss$
by $f_s\colon Y_s\to X_s$
is dense.
Then, $Y\to X$ is \'etale.
\end{lm}

\proof{
If $S$ is regular, 
the assumption that $Y\times_XW\to W$ is \'etale
and Corollary \ref{corSk}.1
implies that the quasi-finite 
morphism  $Y\to X$ of normal 
noetherian schemes
is \'etale in codimension $\leqq 1$.
Since $X$ is regular, the assertion
follows from the purity theorem
of Zariski-Nagata.

Since $X$ and $Y$ are flat over
$S$,
it suffices to show that
for every point $s\in S$,
the morphism
$Y_s=Y\times_Ss
\to X_s$ is \'etale.
Let $S'\to S$ be
the normalization of
the blow-up at the closure of $s\in S$.
Then, there exists a point
$s'\in S'$ above $s\in S$
such that the local ring
${\cal O}_{S',s'}$ is a discrete valuation
ring.
Since the assumption is preserved
by the base change $
{\rm Spec}\, {\cal O}_{S',s'}\to S$,
the morphism
$Y_{s'}=Y\times_S{s'}
\to X_{s'}=X\times_S{s'}$ is \'etale.
Hence 
$Y_s\to X_s$ is also \'etale
as required.
\qed }

\medskip

The following statement is
a combination of the reduced fiber theorem
and the flattening theorem.

\begin{thm}[{\rm \cite[Theorem 2.1$'$]{BLR},
\cite[Th\'eor\`eme (5.2.2)]{RG}}]
\label{thmRFT}
Let $S$ be a noetherian scheme
and $U\subset S$ be
a schematically dense open subscheme.
Let $X$ be a scheme of finite type
over $S$ such that
$X_U=X\times_SU$ is schematically dense
in $X$ and that $X_U\to U$
is flat and reduced.
Then there exists a commutative diagram
\begin{equation}
\begin{CD}
X@<<<X'\\
@VVV@VVV\\
S@<<<S'
\end{CD}
\label{eqRFT}
\end{equation}
of schemes satisfying the following conditions:

{\rm (i)}
The morphism $S'\to S$
is the composition
of a blow-up $S^*\to S$ with center
supported in $S\sm U$
and a faithfully flat morphism
$S'\to S^*$ of finite type
such that
$U'=S'\times_SU\to U$ is \'etale.

{\rm (ii)}
The morphism
$X'\to S'$ is {\em flat} and {\em reduced}.
The induced morphism
$X'\to X\times_SS'$ is finite
and its restriction
$X'\times_{S'}U'\to X\times_SU'$ is an isomorphism.
\end{thm}

If $X_U\to U$ is smooth and if
$S'$ is normal, then
$X'$ is the normalization of $X\times_SS'$
by Corollary \ref{corSk}.1.

For the morphism $S'\to S$
satisfying the condition (i)
in Theorem \ref{thmRFT},
we have a following variant
of the valuative criterion.

\begin{lm}\label{lmvalf}
Let $S$ be a scheme
and $U$ be a dense open subscheme.
Let $S_1\to S$ be a proper morphism
such that $U_1=U\times_SS_1\to U$
is an isomorphism and
let $S'\to S$ be a quasi-finite faithfully
flat morphism.
Let $t\in U$, let
$A\subset K=k(t)$ be a valuation ring
and $T={\rm Spec}\, A\to S$
be a morphism extending $t\to U$.
Then, there exist $t'\in U'=U\times_SS'$
above $t$,
a valuation ring $A'\subset K'=k(t')$
such that $A=A'\cap K$
and a commutative diagram
\begin{equation}
\begin{CD}
T'@>>>S'\\
@VVV@VVV\\
T@>>>S
\end{CD}
\label{eqTS}
\end{equation}
for $T'={\rm Spec}\, A'$.
Further, if
$t=T\times_SU$, then we have
$t'=T'\times_{S'}U'$.
\end{lm}

\proof{
Since $S_1\to S$ is proper
and $U_1\to U$ is an isomorphism,
the morphism
$T\to S$ is uniquely lifted to
$T\to S_1$
by the valuative criterion of
properness.
Let $x_1\in T\times_{S_1}S'$
be a closed point and let
$t'\in t\times_{S_1}S'$
be a point above $t$ such that
$x_1$ is contained in the closure
$T_1=\overline{\{t'\}}
\subset T\times_{S_1}S'$ 
with the reduced scheme structure.
Let $A'\subset k(t')$ be a valuation ring
dominating the local ring ${\cal O}_{T_1,x_1}$.
Then, we have a commutative diagram
(\ref{eqTS}) for $T'={\rm Spec}\, A'$.

Since $t'$ is the unique point of
$t\times_TT'$,
the equality 
$t=T\times_SU$
implies
$t'=T'\times_{S'}U'$.
\qed }

\subsection{Semi-stable curves}

Let $S$ be a scheme.
Recall that a flat separated
scheme $X$ 
of finite presentation
over $S$ is a semi-stable
curve, if every geometric fiber
is purely of dimension $1$
and has at most nodes as
singularities.

\begin{ex}\label{exCD}
{\rm
Let $S$ be a scheme and
$D\subset S$ be an effective Cartier
divisor.
Let $C'\to {\mathbf A}^1_S$
be the blow-up at
$D\subset S
\subset {\mathbf A}^1_S$
regarded as a closed subscheme
by the $0$-section.
Then, the complement
$C_D\subset C'$
of the proper transform
of the $0$-section is
a semi-stable curve over $S$
and is smooth over
the complement $U=S\sm D$.
The exceptional divisor
$\tilde D\subset C_D$
is an effective Cartier divisor 
satisfying
$0\leqq \tilde D\leqq 
D\times_SC_D$.
The difference
$D\times_SC_D-
\tilde D$ equals the proper
transform of ${\mathbf A}^1_D$.

If $S={\rm Spec}\, A$
is affine,
${\mathbf A}^1_S=
{\rm Spec}\, A[t]$
and if $D$ is defined
by a non-zero divisor $f\in A$,
we have 
$C_D={\rm Spec}\, A[s,t]/(st-f)$
and $\tilde D\subset C_D$
is defined by $t$.
}
\end{ex}

\begin{lm}\label{lmsst}
Let $S$ be a scheme
and let $U\subset S$ be
a schematically dense open subscheme.
Let $C$ be a separated
flat scheme of finite presentation 
over $S$ such that the base change
$C_U=C\times_SU$ is a smooth curve
over $U$.
Then, the following
conditions are equivalent:

{\rm (1)}
$C$ is a semi-stable curve over $S$.

{\rm (2)}
Etale locally on $C$ and on $S$,
there exist an effective Cartier
divisor $D\subset S$
such that $D\cap U$
is empty 
and an \'etale morphism
$C\to C_D$ over $S$
to the semi-stable curve $C_D$
defined in Example {\rm \ref{exCD}}.
\end{lm}

\proof{
This is a special case of
\cite[Corollaire 1.3.2]{SGA7}.
\qed }

\medskip

Let $S$ be a normal noetherian scheme
and $j\colon U=S\sm D\to S$ be the 
open immersion of the complement
of an effective Cartier divisor $D$.
Let $i\colon D\to S$ be the closed immersion
and let $\pi_D\colon \bar D\to D$
denote the normalization.
Then, the valuations at the generic points
of irreducible components of $D$
define an exact sequence
$0\to {\mathbf G}_{m,S}
\to j_*{\mathbf G}_{m,U}
\to i_*\pi_{D*}{\mathbf Z}_{\bar D}$
of \'etale sheaves on $S$.

Let $f\colon C=C_D\to S$ be the semi-stable curve over $S$
defined in Example \ref{exCD}.
Let $\tilde j\colon U_C=C\times_SU\to C$
denote the open immersion
and let $\tilde i\colon D_C=C\times_SD\to C$
denote the closed immersion.
Let $A\subset C$ be the exceptional
divisor and $B=D_C-A$
be the proper transform of
${\mathbf A}^1_D$.
Let $a\colon A\to C$ and 
$b\colon B\to C$ 
and $e\colon E=A\cap B
\to C$ denote the closed immersions.
Then, the Cartier divisors
$A,B, D_C\subset C$
defines a commutative diagram
\begin{equation}
\begin{CD}
f^*i_*{\mathbf Z}
@>>> a_*{\mathbf Z}\oplus
b_*{\mathbf Z}
\\
@VVV@VVV\\
f^*(j_*{\mathbf G}_{m,U}/
{\mathbf G}_{m,S})
@>>>
\tilde j_*{\mathbf G}_{m,U_C}/
{\mathbf G}_{m,C}
\end{CD}
\label{eqjAB}
\end{equation}
of \'etale sheaves on $C$.

\begin{pr}\label{prDiv}
Let $S$ be a normal noetherian scheme
and $D\subset S$ be an effective Cartier divisor.
Let $f\colon C=C_D\to S$ be the semi-stable curve
defined in Example {\rm \ref{exCD}}.
Then, the diagram {\rm (\ref{eqjAB})}
induces an exact sequence
\begin{equation}
\begin{CD}
0@>>>
f^*i_*{\mathbf Z}
@>>> 
f^*(j_*{\mathbf G}_{m,U}/
{\mathbf G}_{m,S})
\oplus 
(a_*{\mathbf Z}\oplus
b_*{\mathbf Z})
@>>>
\tilde j_*{\mathbf G}_{m,U_C}/
{\mathbf G}_{m,C}@>>>0
\end{CD}
\label{eqjZ}
\end{equation}
of \'etale sheaves on $D_C$.
\end{pr}

\proof{
Let $z$ be a geometric point of $C$ and we
show the exactness of the stalks of (\ref{eqjZ}) at $z$.
Replacing $S$ by the strict localization
at the image $x$ of $z$,
we may assume that $S$ is strict local and
that $x$ is the closed point.
For $t\in S=S_{(x)}$,
the Milnor fiber $C_{(z)}\times_St$
at $t$ of
the strict localization $C_{(z)}$ at $z$ is geometrically connected
by \cite[Th\'eor\`eme (18.9.7)]{EGA4}.
Further, if $z\in E$ and if $t\in D$,
the fiber at $t$ of $C_{(z)}
\sm E_{(z)}$ has 2
geometrically connected components.

First, we consider the case where $C$ is smooth 
over $S$ at $z$.
Then, since the Milnor fiber $C_{(z),t}$
is connected,
the canonical morphism
$f^*i_*{\mathbf Z}_{\bar D}
\to
i_{C*}{\mathbf Z}_{\bar D_C}$
is an isomorphism.
Hence the stalk of the lower horizontal arrow (\ref{eqjAB}) at $z$
is an injection.
Further this is a surjection by flat descent.

We assume that $C\to S$ is not smooth at $z$.
Let $\tilde D$ be a Cartier divisor of 
$C_{(z)}$ supported on $D_{C_{(z)}}
=C_{(z)}\times_SD$.
Then similarly as above, 
there exists a
Cartier divisor $D_1$ on $S$ supported on $D$
such that $D_0=\tilde D-f^*D_1$ is supported on the
inverse image of $A$.
Define a ${\mathbf Z}$-valued function $n$
on $y\in E_{(z)}=D$ 
as the intersection number
of $D_0$ with the 
fiber $B\times_Sy$.
We show that
the function $n$ is constant.
By adding some multiple of
$A$ to $\tilde D$ if necessary,
we may assume that
$D_0$ is an effective
Cartier divisor of $C$
supported on $A$.
Since $B$ is flat over $D$, 
the pull-back
$D_0\times_CB$ is
an effective Cartier divisor of $B$ finite flat over $D$
by \cite[Proposition (15.1.16) c)$\Rightarrow$b)]{EGA4}.
Hence
the function $n$ is constant. 
Thus 
we have $\tilde D=f^*D_1+n\cdot A$
and the exactness of the stalks of
(\ref{eqjZ}) at $z$ follows.
\qed }

\begin{cor}\label{corDiv}
Let $S$ be a normal noetherian scheme
and $C\to S$ be a semi-stable curve.
Let $x\in S$ be a point
and $z\in C\times_Sx$ be a singular point
of the fiber.
Assume that $z$ is contained in
the intersection of
two irreducible components
$C_1$ and $C_2$ of $C\times_Sx$. 
Let $s_1\colon S\to C$
and $s_2\colon S\to C$ be sections
meeting with the smooth parts
of $C_1$ and $C_2$ respectively.

Let $U\subset S$ be a dense open subscheme
such that $C_U=C\times_SU$
is smooth over $U$ 
and let $\tilde D\subset C$ be an effective
Cartier divisor such that $\tilde D\cap C_U$
is empty.
Define effective Cartier divisors $D_1=s_1^*\tilde D$
and $D_2=s_2^*\tilde D$ of $S$ as
the pull-back of $\tilde D$.

Then, on a neighborhood of
$x$, we have either
$D_1\leqq D_2$ or
$D_2\leqq D_1$.
Suppose we have
$D_1\leqq D_2$ on a neighborhood of $x$.
Then, we have
$D_1\times_SC
\leqq \tilde D\leqq
D_2\times_SC$ on a neighborhood of $z$.
\end{cor}

\proof{
In the notation of the proof of Proposition
\ref{prDiv},
we have $\tilde D=f^*D_1+nA$
for an integer $n$
on an \'etale neighborhood of $z$.
Hence the assertion follows.
\qed }

\medskip
We recall a combination
of a strong version of the semi-stable reduction theorem
for curves over a general base scheme
with the flattening theorem.

\begin{thm}[{\rm \cite[Theorem 2.3.3]{Tst},
\cite[Th\'eor\`eme (5.2.2)]{RG}}]\label{thmsst}
Let $S$ be a noetherian scheme
and $U\subset S$ be 
a schematically dense
open subscheme.
Let $C\to S$ be a separated
morphism of finite type
such that 
$C\times_SU\to U$ is
a smooth curve and that
$C\times_SU\subset U$ is
schematically dense. 
Then, there exists a commutative diagram
$$\begin{CD}
C@<<< C'\\
@VVV@VVV\\
S@<<< S'
\end{CD}$$
of schemes
satisfying the following conditions:

{\rm (i)}
The morphism $S'\to S$ is the composition
of a proper modification $S_1\to S$
such that $U_1=U\times_SS_1\to U$
is an isomorphism
and a faithfully flat morphism 
$S'\to S_1$
such that 
$U'=U\times_SS'\to U_1$ is \'etale
and $U'\subset S'$
is schematically dense.

{\rm (ii)}
The morphism $C'\to S'$ is a 
{\em semi-stable} curve
and
the morphism
$C'\to C\times_SS'$ is 
a proper modification such 
that $C'\times_{S'}U'
\to C\times_SU'$ is 
an isomorphism.
\end{thm}

\begin{cor}\label{cossst}
Let $S$ be a noetherian scheme
and $U\subset S$ be 
a schematically dense open
subscheme.
Let $C\to S$ be a separated morphism
of finite type
such that 
$C_U=C\times_SU\to U$ is a smooth curve
and that
$C_U\subset C$ is schematically dense. 
Let $X\to C$ be a separated morphism
of finite type such that
$X_U=X\times_SU\subset X$
is schematically dense
and that $X_U\to C_U$ is flat and reduced.
Then, there exists a commutative diagram
$$\begin{CD}
X@<<< X'\\
@VVV@VVV\\
C@<<< C'\\
@VVV@VVV\\
S@<<< S'
\end{CD}$$
of schemes
satisfying the following conditions:

{\rm (i)}
The morphism $S'\to S$ is the composition
of a proper modification $S_1\to S$
such that $U_1=U\times_SS_1\to U$
is an isomorphism
and a faithfully flat morphism 
$S'\to S_1$
such that 
$U'=U\times_SS'\to U_1$ is \'etale
and $U'\subset S'$ is schematically dense.

{\rm (ii)}
The morphism $C'\to S'$ is a semi-stable curve
and the morphism
$C'\to C\times_SS'$ is the composition
of a proper modification $C'_0\to C\times_SS'$
such that 
$C'_0\times_{S'}U'\to 
C\times_SU'$ is an isomorphism,
a faithfully flat morphism $C'_1\to C'_0$
such that 
$C'_1\times_{S'}U'\to 
C'_0\times_{S'}U'$
is \'etale
and of a proper modification 
$C'\to C'_1$
such that 
$C'\times_{S'}U'\to 
C'_1\times_{S'}U'$
is an isomorphism.

{\rm (iii)}
The morphism $X'\to C'$
is flat and reduced, 
the morphism $X'\to X\times_CC'$
is finite and 
$X'\times_{S'}U'
\to X\times_CC'\times_{S'}U'$
is an isomorphism.
\end{cor}

\proof{
By the reduced fiber theorem
Theorem \ref{thmRFT}
applied to $X\to C$,
there exists a commutative
diagram 
$$\begin{CD}
X@<<< X_1\\
@VVV@VVV\\
C@<<< C_1
\end{CD}$$
satisfying the conditions
(i) and (ii) loc.\!~cit.
Since $C_1\times_SU\to C\times_SU$
is \'etale
and $C_1\times_SU\subset C_1$
is schematically dense,
by the combination 
Theorem \ref{thmsst} of the
stable reduction theorem
and the flattening theorem,
there exists a commutative
diagram 
$$\begin{CD}
C_1@<<< C_2\\
@VVV@VVV\\
S@<<< S_2
\end{CD}$$
satisfying the conditions 
{\rm (i)} and {\rm (ii)} loc.\! cit.

By the flattening theorem
\cite[Th\'eor\`eme (5.2.2)]{RG}
applied to $S_2\to S$,
there exists a commutative
diagram 
$$\begin{CD}
S_2@<<< S'\\
@VVV@VVV\\
S@<<< S_1
\end{CD}$$
satisfying the condition (i).
We show that
$C'=C_2\times_{S_2}S'$
and 
$X'=X_1\times_{C_1}C'$
satisfy the required conditions.
The base change $C'\to S'$
of a semi-stable curve $C_2\to S_2$
is a semi-stable curve.
Since $C_1\to C$ is obtained
by applying Theorem \ref{thmRFT}
and $C_2\to S_2$ is obtained
by applying Theorem \ref{thmsst},
the composition
$C'=C_2\times_{S_2}S'
\to C'_1=C_1\times_SS'
\to C\times_SS'$
satisfies the condition in (ii).
Finally, the base change
$X'\to C'$ of a flat and reduced morphism
$X_1\to C_1$ is flat and reduced.
Since 
$X'\to C'$ is obtained
by applying Theorem \ref{thmRFT},
the morphism $X'\to X\times_CC'$
satisfies the condition (iii).
\qed
}

\subsection{Subgroups and fiber functor}

For a finite group $G$,
let $(\text{Finite $G$-sets})$ denote
the category of finite sets
with left $G$-action.

\begin{df}\label{dfGC}
We say that a category $C$ 
is a {\em finite Galois category}
if there exist a finite group $G$
and an equivalence
of categories
$F\colon C\to (\text{\rm Finite $G$-sets})$.
If
$F\colon C\to (\text{\rm Finite $G$-sets})$
is an equivalence
of categories, we say that
$G$ is the {\em Galois group}
of the finite Galois category $C$
and call the 
functor $F$ itself or the
composition
$C\to (\text{\rm Finite sets})$
with the forgetful functor
also denoted by
$F$ a {\em fiber functor} of $C$.
\end{df}

We say that a morphism
$F\to F'$ of functors
$F,F'\colon C\to (\text{\rm Finite sets})$
is a surjection if
$F(X)\to F'(X)$
is a surjection for
every object $X$ of $C$.
For a subgroup
$H\subset G$
and for a fiber functor
$F\colon C\to 
(\text{\rm Finite $G$-sets})$,
let $F_H$ denote the functor
$C\to (\text{\rm Finite sets})$
defined by 
$F_H(X)=H\backslash F(X)$.
The canonical morphism
$F\to F_H$ is a surjection.

Surjections $F\to F_H$ are characterized as follows.

\begin{pr}[{\rm cf.\ \cite[Proposition 2.1]{AS}}]\label{prGN}
Let $C$ be a finite Galois category
of Galois group $G$ and
$F\colon C\to (\text{\rm Finite sets})$
be a fiber functor.
Let $F'\colon C\to (\text{\rm Finite sets})$
be another functor
and $F\to F'$
be a surjection of functors.
Then, the following conditions
are equivalent:

{\rm (1)}
For every surjection
$X\to Y$ in $C$,
the diagram 
\begin{equation}
\begin{CD}
F(X)@>>>F'(X)\\
@VVV@VVV\\
F(Y)@>>>F'(Y)
\end{CD}
\label{eqFX}
\end{equation}
is a cocartesian diagram of finite sets.
For every pair of objects
$X$ and $Y$ of $C$,
the morphism
$F'(X)\amalg F'(Y)
\to F'(X\amalg Y)$
is a bijection.

{\rm (2)}
There exists a subgroup
$H\subset G$ such that $F\to F'$
induces an isomorphism
$F_H\to F'$.
\end{pr}

\proof{
(1)$\Rightarrow$(2):
We may assume
$C= (\text{Finite $G$-sets})$
and $F$ is the forgetful functor.
For $X=G$,
the mapping
$F(G)=G\to F'(G)$
is a surjection of finite sets.
Define an equivalence relation
$\sim$ on $G$ by requiring
that $G/\!\!\sim\, \to F'(G)$ to
be a bijection and
set $H=\{x\in G\mid x\sim e\}$.
Then, since the group $G$ acts
on the object $G$ of $C$
by the right action,
the relation $x\sim y$ is equivalent
to $xy^{-1}\in H$.
Since $\sim$ is an equivalence relation,
the transitivity implies
that $H$ is stable under 
the multiplication,
the reflexivity
implies $e\in H$
and the symmetry
implies that $H$
is stable under the inverse.
Hence $H$ is a subgroup
and the surjection
$F(G)=G\to F'(G)$
induces a bijection
$H\backslash G\to F'(G)$.

Let $X$ be an object of $C= (\text{Finite $G$-sets})$
and regard $G\times X$ as a $G$-set
by the left action on $G$.
Then, since the functor
$F'$ preserves the disjoint
union,
we have a canonical isomorphism
$F'(G\times X)
\to F'(G)\times X
\to (H\backslash G)\times X$.
Further,
the cocartesian diagram (\ref{eqFX})
for the surjection 
$G\times X\to X$ in $C$
defined by the action of $G$ is
given by
\begin{equation}
\begin{CD}
G\times X@>>>(H\backslash G)\times X\\
@VVV@VVV\\
X@>>>F'(X).
\end{CD}
\label{eqFG}
\end{equation}
Thus we obtain a bijection
$H\backslash X\to F'(X)$.

The other implication
(2)$\Rightarrow$(1) is clear.
\qed }

\begin{cor}\label{corqt}
Let the notation be as in Proposition {\rm \ref{prGN}}
and let $G'$ be a quotient group.
Let $C'\subset C$ be the full subcategory
consisting of objects $X$ such
that $F(X)$ are $G'$-sets.
Then the subgroup $H'\subset G'$
defined by the surjection 
$F|_{C'}\to F'|_{C'}$ of
the restrictions of the functors
equals the image of $H\subset G$
in $G'$.
\end{cor}

\proof{
If a $G$-set $X$ is a $G'$-set,
the quotient $H\backslash X$
is $H'\backslash X$.
\qed
}

\begin{cor}\label{corGN}
Let $C$ be a finite Galois category
of Galois group $G$
and $F\colon C\to {\text{\rm(Finite $G$-sets)}}$
be a fiber functor.
Let $G'\to G$ be a morphism of groups
and let $F$ 
also denote the functor
$C\to {\text{\rm(Finite $G'$-sets)}}$
defined as the composition
defined by $G'\to G$.
Let $F'\colon C\to {\text{\rm(Finite $G'$-sets)}}$
be another functor
and $F\to F'$ be a surjection of functors
such that the composition
with the forgetful functor
satisfies the condition {\rm (1)}
in Proposition {\rm \ref{prGN}}.

Let $H\subset G$
be the subgroup satisfying
the condition {\rm (2)}
and $G'_1\subset G$ be the image
of $G'\to G$.
Then, the functor $F'$
induces a functor
$C\to {\text{\rm(Finite $G'_1$-sets)}}$
and $G'_1\subset G$ is a subgroup
of the normalizer $N_G(H)$ of $H$.
\end{cor}

\proof{
For an object $X$ of $C$,
$F(X)$ regarded as a $G'$-set
is a $G'_1$-set.
Since $F(X)\to F'(X)$ is a surjection
of $G'$-sets,
$F'(X)$ is also a $G'_1$-set.
Since the left action of $G'_1\subset G$
on the $G$-set
$F(G)=G$ induces an action 
on $F'(G)=H\backslash G$,
the subgroup $H$ is normalized by $G'_1$.
\qed 
}

\section{Dilatations}

\subsection{Functoriality of dilatations}

Let $X$ be a noetherian
scheme and we consider
morphisms
\begin{equation}
\begin{CD}
D@>>> X@<<<Q@<<< Y
\end{CD}
\label{eqDQY}
\end{equation}
of separated schemes of finite 
type over $X$ 
satisfying the following condition:

(i) 
$D\subset X$,
$D_Y=D\times_XY\subset Y$ and
$D_Q=D\times_XQ\subset Q$
are effective Cartier divisors
and $Y\to Q$ is a closed immersion.

\noindent
In later subsections,
we will further assume
the following condition:

(ii) $X$ is normal and
$Q$ is smooth over $X$.

We give examples of 
constructions of $Q$
for a given $Y$ over $X$.

\begin{ex}\label{exQ}
{\rm Assume that $X$ and
$Y$ are separated schemes of finite type
over a noetherian
scheme $S$.

{\rm 1.}
Assume $S={\rm Spec}\, A$ 
and $Y={\rm Spec}\, B$ are
affine.
Then, taking a surjection
$A[T_1,\ldots,T_n]\to B$,
we obtain
a closed immersion
$Y\to Q={\mathbf A}^n_S\times_SX$.

{\rm 2.}
Assume that $Y$ is smooth over $S$.
Then, 
$Q=Y\times_SX\to X$ is smooth
and the canonical morphism
$Y\to Q=Y\times_SX$
is a closed immersion.

{\rm 3.}
Assume that
$\pi\colon Y\to X$ is finite flat
and define a vector bundle
$Q$ over $X$
by the symmetric ${\cal O}_X$-algebra
$S^\bullet \pi_*{\cal O}_Y$.
Then the canonical surjection
$S^\bullet \pi_*{\cal O}_Y
\to \pi_*{\cal O}_Y$
defines a closed immersion
$Y\to Q$.}
\end{ex}

For morphisms
(\ref{eqDQY}) satisfying the
condition (i) above,
we construct a commutative
diagram 
\begin{equation}
\xymatrix{
\bar Y\ar[d]\ar[r]&
Q^{(D)}\ar[d]\ar[dr]&
\\
Y\ar[r]&
Q^{[D]}\ar[r]&
Q}
\label{eqYQD}
\end{equation}
of schemes over $X$ as follows.
Let ${\cal I}_D\subset {\cal O}_X$
and ${\cal I}_Y\subset {\cal O}_Q$
be the ideal sheaves
defining the closed subschemes
$D\subset X$ and $Y\subset Q$.
Let $Q'\to Q$ be the 
blow-up at
$D_Y=D\times_XY
\subset Q$
and define the dilatation
$Q^{[D]}$
at $Y\to Q$ and $D$
to be the largest open subset
of $Q'$ where
${\cal I}_D{\cal O}_{Q'}\supset
{\cal I}_Y{\cal O}_{Q'}$.
Since $D_Y$ is a divisor of $Y$,
by the functoriality of blow-up,
the immersion
$Y\to Q$ is uniquely lifted to a
closed immersion $Y\to Q^{[D]}$.
Let $\bar Y$ and $Q^{(D)}$
be the normalizations
of $Y$ and $Q^{[D]}$
and let $\bar Y\to Q^{(D)}$
be the morphism
induced by the morphism 
$Y\to Q^{[D]}$.
If there is a risk of confusion,
we let $Q^{[D]}$
and $Q^{(D)}$ also be denoted by
$Q^{[D.Y]}$ and $Q^{(D.Y)}$
to make $Y$ explicit.

Locally, 
if $Q={\rm Spec}\, A$
and $Y={\rm Spec}\, A/I$
are affine and if
$D\subset X$ is defined by 
a non-zero divisor $f$,
we have
\begin{equation}
Q^{[D]}={\rm Spec}\, A[I/f]
\label{eqQDAf}
\end{equation}
for the subring $A[I/f]\subset
A[1/f]$ and
the immersion
$Y\to Q^{[D]}$
is defined by the isomorphism
$A[I/f]/(I/f)A[I/f]\to A/I$.

\begin{ex}\label{exQD}
{\rm
Let $X$ be a noetherian scheme
and $D\subset X$ be
an effective Cartier divisor.

{\rm 1.}
Let $Q$ be a smooth separated scheme
over $X$ and $s\colon X\to Q$
be a section.
Let $Y=s(X)\subset Q$
be the closed subscheme.
Then, $Q^{[D]}$ is smooth over $X$.
If $X$ is normal, the canonical morphism
$Q^{(D)}\to Q^{[D]}$
is an isomorphism.

{\rm 2.}
Assume that $X$ is normal.
Let $Q$ be a smooth curve over $X$
and
let $s_1,\ldots,s_n\colon X\to Q$
be sections.
Define
a closed subscheme $Y\subset Q$
as the sum $\sum_{i=1}^ns_i(X)$
of the sections regarded
as effective Cartier divisors of $Q$.
Assume that
$D\subset s_n^*(s_i(X))$
for $i=1,\ldots,n-1$.
Then $Q^{(nD)}\to X$ is smooth and
$Y\times_{Q^{[nD]}}Q^{(nD)}
\subset Q^{(nD)}$
is the sum $\sum_{i=1}^n\tilde s_i(X)$
of the sections $\tilde s_i\colon X\to Q^{(nD)}$
lifting $s_i\colon X\to Q$.

In fact,
we may assume that 
$X={\rm Spec}\, A$ is affine
and,  locally on $Q$, take an \'etale morphism
$Q\to {\mathbf A}^1_X$.
Then, we may assume that
$Q={\mathbf A}^1_X={\rm Spec}\, A[T]$
and $Y$ is defined by
$P=\prod_{i=1}^n(T-a_i)$ for
$a_i\in A$.
We may further assume that
$D$ is defined by a non-zero divisor
$a\in A$ dividing
$a_1,\ldots,a_n$.
Then, we have
$Q^{[nD]}={\rm Spec}\, A[T][P/a^n]$
and 
$T'=T/a$ satisfies
$\prod_{i=1}^n(T'-a_i/a)=P/a^n$
in $A[T][1/a]$.
Hence we have
$Q^{(nD)}={\rm Spec}\, A[T']$
and this equals $Q^{[D.s_n(X)]}$ and is smooth over $X$.
The section $Y\to Q^{[nD]}$
is defined by $P/a^n=0$
and hence
$Y\times_{Q^{[nD]}}Q^{(nD)}
\subset Q^{(nD)}$
is defined by
$A[T']/\prod_{i=1}^n(T'-a_i/a)$.
}
\end{ex}

We study the base change
$Q^{[D]}\times_XD$.

\begin{lm}\label{lmDY}
{\rm 1.}
The canonical morphism
$Q^{[D]}\to Q$ induces
\begin{equation}
Q^{[D]}\times_XD
=
Q^{[D]}\times_QD_Y
\to D_Y.
\label{eqQDDY}
\end{equation}

{\rm 2.}
If $Y\to Q$ is a regular immersion
and if $T_YQ$ and
$T_DX$ denote the normal bundles,
we have a canonical isomorphism
\begin{equation}
T_YQ(-D_Y)\times_YD_Y
=
(T_YQ\times_YD_Y)
\otimes
(T_DX\times_DD_Y)^{\otimes -1}
\to
Q^{[D]}\times_XD.
\label{eqTDY}
\end{equation}
The isomorphism {\rm (\ref{eqTDY})}
depends only on
the restriction $D_Y\to Q$
and not on $Y\to Q$
itself.

{\rm 3.}
Assume that
$Q$ is smooth over $X$
and 
$X=Y\to Q$ is a section.
Let $T(Q/X)$ denote the 
relative tangent bundle
defined by the symmetric 
${\cal O}_Q$-algebra
$S_{{\cal O}_D}^\bullet
\Omega^{1}_{Q/X}$.
Then,
we have a canonical isomorphism
\begin{equation}
T(Q/X)(-D)\times_QD
=
(T(Q/X)\times_QD)
\otimes
T_DX^{\otimes -1}
\to
Q^{[D]}\times_XD.
\label{eqTDX}
\end{equation}
The isomorphism {\rm (\ref{eqTDX})}
depends only on
the restriction $D\to Q$
and not on the section $X\to Q$
itself.
\end{lm}

\proof{
1. Since 
${\cal I}_D{\cal O}_{Q^{[D]}}\supset
{\cal I}_Y{\cal O}_{Q^{[D]}}$ on $Q^{[D]}$
by the definition of $Q^{[D]}$,
we have
$Q^{[D]}\times_XD
=
Q^{[D]}\times_QD_Y$.
Hence, we
obtain a morphism
$Q^{[D]}\times_XD\to D_Y$.

2. Assume that $Y\to Q$
is a regular immersion.
Then, $D_Y\to Q$ is also
a regular immersion and 
the normal bundle $T_{D_Y}Q$
fits in an exact sequence
$0\to T_{D_Y}D_Q\to
T_{D_Y}Q
\to T_DX\times_DD_Y\to 0$
depending only on 
$D\to X$ and $D_Y\to Q$
and not on $Y\to Q$.
Let $Q'\to Q$ be the blow-up 
at $D_Y\subset Q$.
Then, the exceptional divisor
$Q'\times_QD_Y$ is canonically
identified with the projective
space bundle 
${\mathbf P}(T_{D_Y}Q)$
over $D_Y$.
Its open subset
$Q^{[D]}\times_QD_Y$
is identified as in (\ref{eqTDY})
since $T_{D_Y}D_Q
=T_YQ\times_YD_Y$.

3.
Since the normal bundle
$T_XQ$ is canonically identified
with the restriction
$T(Q/X)\times_QX$
of the relative tangent bundle,
the assertion follows from 2.
\qed
}
\medskip

We give a sufficient condition
for the morphism
$\bar Y\to Q^{(D)}$
to be an immersion.

\begin{lm}\label{lmb+}
Assume that $X$ and $Y\sm D_Y$
are normal and 
let $\pi\colon \bar Y\to Y$
be the normalization.
Assume that 
$\bar Y\to X$ is \'etale and that
$\pi_*{\cal O}_{\bar Y}/{\cal O}_Y$
is an ${\cal O}_{D_Y}$-module.
Then,  the finite morphism
$\bar Y\to Q^{(2D)}$ is a closed immersion.
\end{lm}

\proof{
Since the assertion is \'etale local on $Y$,
we may assume that
$Y\to X$ is finite and that
the \'etale covering $\bar Y\to X$
is split.
We may further assume that $X,Y$ and $Q$ 
are affine and that $D$
is defined by a non-zero divisor $f$ on $X$. 
Let $Y={\rm Spec}\, A$,
$\bar Y={\rm Spec}\, \bar A$,
$Q={\rm Spec}\, B$,
$Q^{[2D]}={\rm Spec}\, B^{[2D]}$,
$Q^{(2D)}={\rm Spec}\, B^{(2D)}$
for $A=B/I, B^{[2D]}=
B[I/f^2]\subset B[1/f]$ 
and the normalization $B^{(2D)}$
of $B^{[2D]}$.
Since $\bar Y\to X$ is a split \'etale covering,
it suffices to show that for
every idempotent $e\in \bar A$,
there exists a lifting $\tilde e\in B^{(2D)}$.

Since $\bar A/A$ is annihilated by $f$,
the product $f e=g$ is an element of $A$.
Let $\tilde g\in B$ be a lifting of $g$.
Since $e^2=e$,
the element
$h=\tilde g^2-f\tilde g\in B$
is contained in $I$
and hence $h/f^2\in B[1/f]$ is an element of
$B^{[2D]}$.
Thus $\tilde e=\tilde g/f \in B[1/f]$
is a root of the polynomial 
$T^2-T-h/f^2\in B^{[2D]}[T]$
and is an element of $B^{(2D)}$.
Since $\tilde e$ is a lifting of $e$,
the assertion follows.
\qed 
}

\medskip

We study the functoriality of the construction.
We consider
a commutative diagram
\begin{equation}
\begin{CD}
D\times_XX'&\, 
\subset D'\subset\,
& X'@<<<Q'@<<< Y'\\
@VVV@VVV@VVV@VVV\\
D&\subset& X@<<<Q@<<< Y
\end{CD}
\label{eqFun}
\end{equation}
of schemes such that the both lines
satisfy the condition (i)
on the diagram (\ref{eqDQY}).
Then, by the functoriality of dilatations
and normalizations,
we obtain a commutative diagram
\begin{equation}
\begin{CD}
Y'@>>>
Q^{\prime[D']}@<<< Q^{\prime(D')}@<<<\bar Y'\\
@VVV@VVV@VVV@VVV\\
Y@>>>
Q^{[D]}@<<< Q^{(D)}@<<<\bar Y.
\end{CD}
\label{eqFun[]}
\end{equation}
The diagram {\rm (\ref{eqFun[]})} 
induces a morphism
\begin{equation}
Q^{\prime (D')}
\times_{Q^{\prime [D']}}Y'
\to
Q^{(D)}
\times_{Q^{[D]}}Y.
\label{eqQDY}
\end{equation}
Let $\bar x$ be a geometric point
of $D$ and 
$\bar x'$ be a geometric point
of $D\times_XX'$ above $\bar x$.
Then the diagram {\rm (\ref{eqFun[]})} 
also induces a mapping
\begin{equation}
\pi_0(Q^{\prime (D')}_{\ \bar x'})
\to
\pi_0(Q^{(D)}_{\ \bar x})
\label{eqpiDY}
\end{equation}
of the sets of connected components
of the geometric fibers.

First we study the dependence on $Q$.

\begin{pr}\label{prFunQ}
Suppose $X=X',Y=Y'$ and $D=D'$
and let $\bar x$ be a geometric point of $D$.

{\rm 1.}
Assume that $Q$ and $Q'$
are smooth over $X$.
Then, the square
\begin{equation}
\begin{CD}
Q^{\prime[D]}@<<<Q^{\prime(D)}\\
@VVV@VVV\\
Q^{[D]}@<<< Q^{(D)}
\end{CD}
\label{eqQD}
\end{equation}
is cartesian.
The induced morphism
$Q^{\prime (D)}
\times_{Q^{\prime [D]}}Y
\to
Q^{(D)}
\times_{Q^{[D]}}Y$
{\rm (\ref{eqQDY})}
is an isomorphism over $Y$
and the induced mapping
$\pi_0(Q^{\prime(D)}_{\ \bar x})
\to 
\pi_0(Q^{(D)}_{\ \bar x})$
{\rm (\ref{eqpiDY})} is a bijection.

{\rm 2.}
Assume that $Q'\to Q$ is smooth
and let $T=T(Q'/Q)$ denote the relative
tangent bundle of $Q'$ over $Q$.
Then $Q^{\prime[D]}\to Q^{[D]}$ 
is also smooth and 
there exists a cartesian diagram
\begin{equation}
\begin{CD}
T(-D)\times_{Q'}D_Y
@<<<
Q^{\prime[D]}\times_XD\\
@VVV
\hspace{-20mm}
\square
\hspace{17mm}
@VVV\\
D_Y
@<{\rm (\ref{eqQDDY})}<<
Q^{[D]}\times_XD.
\end{CD}
\label{eqQDD}
\end{equation}
\end{pr}

\proof{
2.
First, we show the case where 
$Q'\to Q$ admits a section $Q\to Q'$
extending $Y\to Q'$.
The section $Q\to Q'$ defines 
a section $Q^{[D]}\to Q'\times_QQ^{[D]}$.
Define $(Q'\times_QQ^{[D]})^{[D_{Q^{[D]}}.Q^{[D]}]}$
to be the dilatation of
$Q'\times_QQ^{[D]}$ for the section
$Q^{[D]}\to Q'\times_QQ^{[D]}$
and a divisor $D_{Q^{[D]}}=D\times_XQ^{[D]}$ 
over $Q^{[D]}$.
We show that the canonical morphism
$Q^{\prime [D]}\to Q'\times_QQ^{[D]}$
induces an isomorphism
\begin{equation}
Q^{\prime [D]}\to 
(Q'\times_QQ^{[D]})^{[D_{Q^{[D]}}.Q^{[D]}]}.
\label{eqQQ'D}
\end{equation}

Since the question is \'etale local on $Q'$,
we may assume that 
$Q'={\mathbf A}^n_Q$
and the section $Q\to Q'$
is the $0$-section.
Further, we may assume that
$Q={\rm Spec}\, A$ and
$Y={\rm Spec}\, A/I$
are affine and that $D\subset X$
is defined by a non-zero divisor $f$ on $X$.
We set
$A'=A[T_1,\ldots,T_n]$
and $Q'={\rm Spec}\, A'$.
The $0$-section $Q\to Q'$
is defined by the ideal
$J=(T_1,\ldots,T_n)\subset A'$.
We have
$Q^{[D]}=
{\rm Spec}\, A[I/f]$
and 
$Q^{\prime [D]}=
{\rm Spec}\, A'[I'/f]$
for $I'=IA'+J$.
Since 
$A'[I'/f]=A[I/f][T_1/f,\ldots,T_n/f]$
as a subring of $A'[1/f]$,
we obtain an isomorphism
(\ref{eqQQ'D}).

By the isomorphism (\ref{eqQQ'D})
and Example \ref{exQD}.1,
the morphism $Q^{\prime [D]}\to Q^{[D]}$
is smooth.
Further by Lemma \ref{lmDY}.3,
we obtain a cartesian
diagram (\ref{eqQDD}),
depending only on $D\to X$,
$D_Y\to Q$ and
$D_Y\to Q'$ but 
not on the choice of section
$Q\to Q'$ extending $Y\to Q'$.

We prove the general case.
Since $Q'\to Q$ has a section
on $Y\subset Q$,
locally on $Q$,
there exist a closed
subscheme $Q_1\subset Q'$ 
\'etale over $Q$
such that $Y\to Q'$ is 
induced by $Y\to Q_1$.
For the smoothness 
of $Q^{\prime [D]} \to Q^{[D]}$,
since the assertion is \'etale local,
we may assume that $Q_1=Q$ is a section.
Hence the smoothness
$Q^{\prime [D]} \to Q^{[D]}$ follows.
Further since
the cartesian diagram (\ref{eqQDD})
defined \'etale locally
is independent of the choice of section,
we obtain (\ref{eqQDD}) for $Q'$
by patching.

1.
First, we show the case where
$Q'\to Q$ is smooth.
Then by 2, 
$Q^{\prime [D]}\to Q^{[D]}$ is 
also smooth and
the fibered product
$Q^{(D)}\times_{Q^{[D]}}Q^{\prime [D]}$ is normal.
Hence the square (\ref{eqQD}) is cartesian
and
the morphism (\ref{eqQDY}) is an isomorphism.
By the cartesian squares
(\ref{eqQD})
and (\ref{eqQDD}),
$Q^{\prime (D)}_{\ \bar x}$
is a vector bundle over
$Q^{(D)}_{\ \bar x}$.
Hence (\ref{eqpiDY})
is a bijection.

We show the general case.
A morphism $f\colon Q'\to Q$ is decomposed
as the composition of the projection
${\rm pr}_2\colon
Q'\times_XQ\to Q$
and a section of the projection
${\rm pr}_1\colon
Q'\times_XQ\to Q'$.
Hence, the cartesian squares (\ref{eqQD})
and the bijections
(\ref{eqpiDY})
for the projections imply
those for $f$
respectively.
The cartesian square (\ref{eqQD})
for $f$ implies an isomorphism
(\ref{eqQDY}) for $f$.
\qed }

\medskip

\begin{cor}\label{corFunQ}
Assume that $Q$ and $Q'$
are smooth over $X$.
Then, the morphism
$Q^{\prime (D')}
\times_{Q^{\prime [D']}}Y'
\to
Q^{(D)}
\times_{Q^{[D]}}Y$
{\rm (\ref{eqQDY})} 
is independent of $Q'\to Q$.
Let $\bar x$ be a geometric point
of $D$ and 
$\bar x'$ be a geometric point
of $D'$ above $\bar x$.
Then the mapping
$\pi_0(Q^{\prime (D')}_{\ \bar x'})
\to
\pi_0(Q^{(D)}_{\ \bar x})$
{\rm (\ref{eqpiDY})} 
is independent of
morphism $Q'\to Q$.
\end{cor}

\proof{
Decompose a morphism
$Q'\to Q$ as
$Q'\to Q'\times_XQ\to Q$.
Then the isomorphism {\rm (\ref{eqQDY})} 
and the bijection (\ref{eqpiDY}) for
$Q'\to Q'\times_XQ$
are the inverses of those
for the projection
$Q'\times_XQ\to Q'$.
Hence the assertion follows.
\qed }

\medskip

By the canonical isomorphism
(\ref{eqQDY}),
the finite scheme
$Y\times_{Q^{[D]}}Q^{(D)}$
over $Y$ is independent of $Q$.
We let it denoted by
$Y^{(D)}$.

\begin{lm}\label{lmFunY}
Suppose that the squares
$$
\begin{CD}
D'@>>> X'\\
@VVV
\hspace{-10mm}
\square
\hspace{7mm}
@VVV\\
D@>>>X,
\end{CD}\qquad
\begin{CD}
Q'@<<< Y'\\
@VVV
\hspace{-10mm}
\square
\hspace{7mm}
@VVV\\
Q@<<< Y
\end{CD}$$
are cartesian.

{\rm 1.}
The morphism
$Q^{\prime[D']}\to
Q^{[D]}\times_QQ'$
is a closed immersion and
$Q^{\prime(D')}\to
Q^{(D)}\times_QQ'$
is finite.
Consequently, the morphism
$Q^{\prime(D')}\times_{Q'}Y'
\to
Q^{(D)}\times_QY$
is finite if $Y'\to Y$ is finite.
Further, if $Q$ and $Q'$ are normal,
then $Q^{\prime(D')}$
equals the normalization of
$Q^{(D)}\times_QQ'$
in $Q'\sm D'\times_{X'}Q'$.

{\rm 2.}
If $Q'\to Q$ is flat,
the square
$$\begin{CD}
Q^{\prime[D']}@>>> Q'\\
@VVV
\hspace{-10mm}
\square
\hspace{7mm}
@VVV\\
Q^{[D]}@>>> Q
\end{CD}$$
is cartesian.
\end{lm}

\proof{
Since the assertion is local on a neighborhood
of $Y'\subset Q'$,
we may assume that
$Q={\rm Spec}\, A$,
$Y={\rm Spec}\, A/I$,
$Q'={\rm Spec}\, A'$ and
$Y'={\rm Spec}\, A'/IA'$
are affine and that
$D$ is defined by a non-zero divisor $f$ on $X$.
Then, we have
$Q^{[D]}={\rm Spec}\, A[I/f]$
and 
$Q^{\prime [D']}={\rm Spec}\, A'[IA'/f]$.

1.
Since
$A[I/f]\otimes_AA'\to
A'[IA'/f]$ is a surjection,
the morphism
$Q^{\prime[D']}\to
Q^{[D]}\times_QQ'$
is a closed immersion.
The remaining assertions follow from this immediately.

2. 
If $A\to A'$ is flat,
the injection $A[I/f]\to A[1/f]$
induces an injection
$A'\otimes_AA[I/f]
\to A'\otimes_AA[1/f]=A'[1/f]$.
Hence the surjection
$A'\otimes_AA[I/f]
\to A'[IA'/f]$ is an isomorphism.
\qed }

\medskip

The construction of $Q^{(D)}$
commutes with base change
if $Q^{(D)}\to X$ is flat and reduced.

\begin{lm}\label{lmFunX}
Suppose that the diagram
{\rm (\ref{eqFun})} is cartesian
and $D'=D\times_XX'$.
Assume that one
of the following conditions
is satisfied:

{\rm (i)}
$X'$ is normal,
$Q\to X$ is smooth
and 
$Q^{(D)}\to X$
is flat and reduced.

{\rm (ii)}
$X'\to X$ is smooth.

\noindent
Then the square
\begin{equation}
\begin{CD}
Q^{(D)}@<<<Q^{\prime(D')}\\
@VVV@VVV\\
X@<<<X'
\end{CD}
\label{eqFunX}
\end{equation}
is cartesian.
\end{lm}

\proof{
By Lemma \ref{lmFunY}.1,
$Q^{\prime (D')}$
is the normalization of
$Q^{(D)}\times_XX'$.
If the condition (i)
is satisfied, then
$Q^{(D)}\times_XX'$
is normal
by Corollary \ref{corSk}.1.
If $X'\to X$ is smooth,
then
$Q^{(D)}\times_XX'$
is smooth over $Q^{(D)}$
and is normal.
Hence the square (\ref{eqFunX})
is cartesian in both cases.
\qed }

\medskip

We study the dependence on $D$
and show that the canonical morphism
contracts the closed fiber.

\begin{lm}\label{lmFunD}
Suppose $X=X',Y=Y'$ and $Q=Q'$
and that $D_1=D'-D$
is an effective Cartier divisor of $X$.
Then, the morphism
$Q^{[D']}\to Q^{[D]}$
(resp.\ $Q^{(D')}\to Q^{(D)}$)
induces a morphism
$Q^{[D']}\times_QD_{1,Y}\to 
D_{1,Y}\subset Y\subset Q^{[D]}$
(resp.\ $Q^{(D')}\times_QD_{1,Y}\to 
Q^{(D)}\times_{Q^{[D]}}D_{1,Y}
\subset Q^{(D)}$).
\end{lm}

\proof{
We consider the immersion 
$Y\to Q^{[D]}$
lifting $Y\to Q$.
Then, the morphism
$Q^{[D']}\to Q^{[D]}$
induces an isomorphism
$Q^{[D']}\to (Q^{[D]})^{[D_1]}$
to the dilatation
$(Q^{[D]})^{[D_1]}$ of
$Q^{[D]}$ for $Y\to Q^{[D]}$
and $D_1\subset X$.
Hence the morphism (\ref{eqQDDY})
defines a morphism
$Q^{[D']}\times_QD_{1,Y}
\to D_{1,Y}$.
The assertion for $Q^{(D')}$
follows from this.
\qed }

\subsection{Dilatations and complete intersection}

We give a condition for the right square
in (\ref{eqFun}) to be cartesian.

\begin{lm}\label{lmlci}
Let $S$ be a noetherian scheme
and let $Q\to P$ be a quasi-finite
morphism of smooth schemes of
finite type over $S$.
If $Q\to P$ is flat
on dense open subschemes,
then $Q\to P$ is flat and 
locally of complete intersection
of relative virtual dimension $0$.
\end{lm}

\proof{
Let $U\subset P$
and $V\subset Q$
be dense open subschemes
such that $V\to U$ is flat.
Then the relative dimension
of $V\to S$ is the same
as that of $U\to S$.
Hence, we may assume
that the relative dimensions of
$P\to S$ and 
$Q\to S$ are the same integer $n$.

The morphism $Q\to P$
is the composition of the graph
$Q\to Q\times_SP$ and the projection
$Q\times_SP\to P$.
For every point $x\in P$,
the fiber
$Q\times_Px\to Q\times_Sx$
is a regular immersion of codimension $n$.
Hence by \cite[Proposition (15.1.16)
c)$\Rightarrow$b)]{EGA4}
applied to the immersion
$Q\to Q\times_SP$ over $P$,
the immersion 
$Q\to Q\times_SP$
is also a regular 
immersion of codimension $n$
and $Q\to P$ is flat.
\qed }

\begin{lm}\label{lmPQ}
Let $S$ be a noetherian scheme
and let $Y\to X$ be a
morphism of schemes of
finite type over $S$.

{\rm 1.}
Suppose that there exists a cartesian diagram
\begin{equation}
\begin{CD}
Q@<<< Y\\
@VVV
\hspace{-10mm}
\square
\hspace{7mm}
@VVV\\
P@<<< X
\end{CD}
\label{eqPQ}
\end{equation}
of schemes of finite type over $S$
satisfying the following conditions:

$P$ and $Q$ are smooth over $S$
and $Q\to P$ is quasi-finite
and is flat
on dense open subschemes.
The horizontal arrows are closed immersions.

\noindent
Then $Y\to X$ is quasi-finite, flat and 
locally of complete intersection
of relative virtual dimension $0$.

{\rm 2.}
Conversely, suppose that
$Y\to X$ is 
finite (resp.\ quasi-finite) and
locally of complete
intersection of relative virtual dimension $0$.
Then $Y\to X$ is flat
and, locally on $X$
(resp. locally on $X$ and on $Y$),
there exists a cartesian diagram
{\rm (\ref{eqPQ})} satisfying the following conditions:

$P$ and $Q$ are smooth 
of the same relative dimension over $S$
and $Q\to P$ is quasi-finite
and flat.
The horizontal arrows are closed immersions.
\end{lm}

\proof{
1.
By Lemma \ref{lmlci},
the quasi-finite morphism
$Q\to P$ is flat and locally of complete intersection.
Hence
$Y\to X$ is also quasi-finite,
flat and locally of complete intersection
of relative virtual dimension $0$.

2.
Since the assertion is local,
we may assume that $S, X$
and $Y$ are affine.
Take a closed immersion
$Q_1={\mathbf A}^m_X\gets Y$.
Since the immersion
$Y\to Q_1$ is a regular immersion
of codimension $m$
and since
$Y\to X$ is finite (resp.\ quasi-finite),
after shrinking $X$ (resp.\ $Q_1$ and $Y$), we may
assume that the ideal defining
$Y\subset Q_1$
is generated by $m$ sections
$f_1,\ldots,f_m$ of ${\cal O}_{Q_1}$.
Also take a closed immersion 
$P_1={\mathbf A}^n_S\gets X$
and an open subscheme 
$Q\subset {\mathbf A}^m_{P_1}$
to obtain a cartesian diagram
\begin{equation}
\begin{CD}
Q@<<< Q_1@<<<Y\\
@VVV
\hspace{-10mm}
\square
\hspace{7mm}
@VVV\\
P_1@<<<X.
\end{CD}
\label{eqcocaP1}
\end{equation}
Taking sections
$\tilde f_1,\ldots,\tilde f_m$ of ${\cal O}_Q$
lifting $f_1,\ldots,f_m$ after shrinking $Q$
if necessary,
define a morphism
$Q\to P={\mathbf A}^m_{P_1}$.
Then, we obtain a cartesian diagram
\begin{equation}
\begin{CD}
Q@<<< Q_1@<<<Y\\
@VVV
\hspace{-10mm}
\square
\hspace{7mm}
@VVV
\hspace{-10mm}
\square
\hspace{7mm}
@VVV\\
P@<<<{\mathbf A}^m_X
@<<<X
\end{CD}
\label{eqPQ1}
\end{equation}
where the lower right horizontal arrow
${\mathbf A}^m_X\to X$
is the $0$-section.

The schemes
$P={\mathbf A}^{n+m}_S$
and $Q\subset {\mathbf A}^{n+m}_S$
are smooth over $S$.
Since $Y\to X$
is quasi-finite,
after replacing $Q$ 
by a neighborhood of
$Y$ if necessary,
the morphism $Q\to P$
is quasi-finite.
Since $Q$ and $P$
are smooth of the same relative dimension
over $S$,
the morphism $Q\to P$
is flat on dense open subschemes.
By Lemma \ref{lmlci},
the quasi-finite morphism
$Q\to P$ is flat and 
hence $Y\to X$ is also flat.
\qed }

\medskip

We give examples of construction
of the diagram (\ref{eqPQ}).

\begin{ex}\label{exPQ}
{\rm 
Assume that $X$ and
$Y$ are schemes of finite type
over a noetherian
scheme $S$.

{\rm 1.}
Assume $X={\rm Spec}\, A$ 
and $Y={\rm Spec}\, B$ are
affine.
Let $A[T_1,\ldots,T_n]
/(f_1,\ldots,f_n)\to B$
be an isomorphism
and define a morphism
$Q={\mathbf A}^n_X
={\rm Spec}\, A[T_1,\ldots,T_n]
\to 
P={\mathbf A}^n_X$
by $f_1,\ldots,f_n$.
Then, we obtain
a cartesian diagram
{\rm (\ref{eqPQ})}
by defining the section
$X\to P={\mathbf A}^n_X$ to be
the $0$-section.

{\rm 2.}
Assume that $X$ and $Y$ 
are smooth over a noetherian scheme
$S$.
Then, we obtain a cartesian
diagram
$$\begin{CD}
Y@>>> Q&\,=Y\times_SX\\
@VVV
\hspace{-10mm}
\square
\hspace{7mm}
@VVV\\
X@>>> P&\,=X\times_SX.
\end{CD}$$}
\end{ex}

\medskip

Assume that
$Q^{(D)}\to P^{(D)}$
is \'etale on a neighborhood of
$Q^{(D)}\times_XD$.
Let $\bar x$ be a geometric point of $D$
and let $0_{\bar x}$ denote
the geometric point above the origin of
the vector space $P^{(D)}_{\ \bar x}$ 
over $\bar x$.
Then, since 
$Q^{(D)}_{\ \bar x}\to P^{(D)}_{\ \bar x}$
is finite \'etale,
we have an action of the fundamental
group $\pi_1(P^{(D)}_{\ \bar x},0_{\bar x})$
on $Y^{(D)}_{\ \bar x}
=Q^{(D)}_{\ \bar x}
\times_{P^{(D)}_{\ \bar x}}0_{\bar x}$.
The action on 
$Y^{(D)}_{\ \bar x}$ is compatible
with the canonical mapping
$Y^{(D)}_{\ \bar x}\to
\pi_0(Q^{(D)}_{\ \bar x})$
with respect to the trivial action on
$\pi_0(Q^{(D)}_{\ \bar x})$
and is transitive on the inverse image
of each element of
$\pi_0(Q^{(D)}_{\ \bar x})$.

Since $Q^{[D]}\to P^{[D]}\times_PQ$
is an isomorphism
by Lemma \ref{lmFunY}.2,
for a geometric point $\bar y$
of $Y_{\bar x}$
and for the geometric point $0_{\bar y}$
of $Q^{[D]}_{\ \bar y}$
above $P^{(D)}_{\ \bar x}$,
we have canonical isomorphisms
$Q^{[D]}_{\ \bar y}
=
Q^{[D]}\times_Q{\bar y}
\to 
P^{[D]}_{\ \bar x}=
P^{(D)}_{\ \bar x}$
and
$\pi_1(Q^{[D]}_{\ \bar y},0_{\bar y})
\to 
\pi_1(P^{(D)}_{\ \bar x},0_{\bar x})$.
The action of
$\pi_1(P^{(D)}_{\ \bar x},0_{\bar x})$
on $Y^{(D)}_{\ \bar x}$
is compatible with the action of
$\pi_1(Q^{[D]}_{\ \bar y},0_{\bar y})$
on $Y^{(D)}_{\ \bar x}\times_{Y_{ \bar x}}
\bar y$.
For a morphism $Q'\to Q$,
the canonical morphism
$\pi_1(Q^{\prime [D]}_{\ \bar y},0_{\bar y})
\to \pi_1(Q^{[D]}_{\ \bar y},0_{\bar y})$
is compatible with the actions
on $Y^{(D)}_{\ \bar x}\times_{Y_{ \bar x}}
\bar y$.

We study the relation between
the \'etaleness of
$Q^{(D)}\to P^{(D)}$
and the annihilator of
${\cal O}_{Y^{(D)}}\otimes
_{{\cal O}_Y}\Omega^1_{Y/X}$.

\begin{lm}\label{lmetale}
Let
$$\begin{CD}
Q@<<< Y\\
@VVV
\hspace{-10mm}
\square
\hspace{7mm}
@VVV\\
P@<<<X
\end{CD}$$
be a cartesian diagram
of separated schemes 
of finite type over $X$.
Assume that
$P$ and $Q$ are smooth over
$X$ and
that the vertical arrows
are quasi-finite and flat.

Assume that there exists
an effective Cartier divisor
$D_1\subset D=D_1+D_0$ 
of $X$ such that
${\cal O}_{Y^{(D)}}
\otimes_{{\cal O}_Y}\Omega^1_{Y/X}$
is annihilated by ${\cal I}_{D_1}
\subset {\cal O}_X$
and that we have an equality
$D_0=D$ of underlying sets.
Then, there exists an
open neighborhood 
$W\subset Q^{[D]}$ of
$Q^{[D]}\times_XD$ such that
$Q^{[D]}\to P^{[D]}$
is \'etale on 
$W\sm(Q^{[D]}\times_XD)$.
\end{lm}

\proof{
It suffices to show that 
each irreducible component
$Z\subset Q^{[D]}$ 
of the inverse image
of the support of 
$\Omega^1_{Q/P}$
is either a subset of
$Q^{[D]}\times_XD$
or does not meet
$Q^{[D]}\times_XD$,
since $Q^{[D]}\to Q$
is an isomorphism on the complement
of the inverse images of $D$.
Assume that $Z$ is not
a subset of $Q^{[D]}\times_XD$
but does meet $Q^{[D]}\times_XD$
and regard $Z$ as an integral
closed subscheme of $Q^{[D]}$.
Then, $D\times_XZ\subset Z$
is a non-empty effective 
Cartier divisor.

Since the assertion is \'etale local on $Y$,
we may assume that
$Y\to X$ is faithfully flat and finite.
Let $T_0\subset Z\times_XY^{(D)}$
be the closure 
of the complement
$Z\times_XY^{(D)}
\sm D\times_X(Z\times_XY^{(D)})$
and $T$ be its normalization.
Then, since $Y\to X$ is finite
surjective,
$T\to Z$ is also finite surjective.
Hence $D_T=D\times_XT\subset T$
is a non-empty effective 
Cartier divisor.

By the assumption that
${\cal O}_{Y^{(D)}}
\otimes_{{\cal O}_Y}\Omega^1_{Y/X}$
is annihilated by ${\cal I}_{D_1}
\subset {\cal O}_X$,
the ${\cal O}_T$-module
${\cal O}_T\otimes_{{\cal O}_Y}
\Omega^1_{Y/X}$
is annihilated by
${\cal I}_{D_1}\cdot {\cal O}_T$.
Since $D_T$ is a scheme over
$Q^{[D]}\times_XD$,
we have an isomorphism
${\cal O}_{D_T}\otimes_{{\cal O}_Q} 
\Omega^1_{Q/P}
\to
{\cal O}_{D_T}\otimes_{{\cal O}_Y}
\Omega^1_{Y/X}$
by Lemma \ref{lmDY}.1.
Thus
${\cal O}_{D_T}\otimes_{{\cal O}_Q} 
\Omega^1_{Q/P}$
is also annihilated by
${\cal I}_{D_1}\cdot {\cal O}_{D_T}$.
Since $D=D_1+D_0$,
this means an inclusion
${\cal I}_{D_1}\cdot
{\cal O}_T\otimes_{{\cal O}_Q} 
\Omega^1_{Q/P}
\subset
{\cal I}_{D_0}\cdot
{\cal I}_{D_1}\cdot
{\cal O}_T\otimes_{{\cal O}_Q} 
\Omega^1_{Q/P}$.
By Nakayama's lemma,
we have
${\cal I}_{D_1}\cdot
{\cal O}_T\otimes_{{\cal O}_Q} 
\Omega^1_{Q/P}=0$
on a neighborhood of
$D_0\times_XT$.

Since $Z$ is a subset of
the inverse image of
support of $\Omega^1_{Q/P}$,
the annihilator ideal of
${\cal O}_T\otimes_{{\cal O}_Q} 
\Omega^1_{Q/P}$ is $0$.
This contradicts to that
$D_0\times_XT=D_T$
is non-empty.
\qed }

\begin{lm}\label{lmet}
Assume $X$ is normal
and let
$$\begin{CD}
Q@<<< Y\\
@VVV
\hspace{-10mm}
\square
\hspace{7mm}
@VVV\\
P@<<<X
\end{CD}$$
be a cartesian diagram
of separated schemes 
of finite type over $X$.
Assume that
$P$ and $Q$ are smooth over
$X$ and
that the vertical arrows
are quasi-finite and flat.

Let $Y_0$
be a closed subscheme
of $Y$ \'etale over $X$ 
satisfying an equality $D_{Y_0}=
D_Y$ of underlying sets
and let ${\cal J}_0
\subset {\cal O}_{D_Y}$
be the nilpotent ideal defining
$D_{Y_0}\subset D_Y$.
Let $n\geqq 1$ be
an integer satisfying
${\cal J}_0^n=0$ and let
$D_0\subset D$
be an effective Cartier divisor
on $X$
satisfying $nD_0\leqq D$.

Assume that $Y^{(D)}
=Y\times_{Q^{[D]}}Q^{(D)}$
is \'etale over $X$.
Then ${\cal O}_{Y^{(D)}}
\otimes_{{\cal O}_Y}
\Omega^1_{Y/X}$
is annihilated by the ideal
${\cal I}_{D-D_0}\subset 
{\cal O}_X$ defining
$D-D_0\subset X$.
\end{lm}

\proof{
Let ${\cal I}\subset
{\cal I}_0\subset {\cal O}_Q$
and ${\cal I}_D
\subset {\cal I}_{D_0}
\subset {\cal O}_X$ 
be the ideals
defining the closed subschemes
$Y_0\subset Y\subset Q$
and $D_0\subset D\subset X$.
Let $Y_0^{(n)}\subset Q$
denote the closed scheme
defined by the ideal
${\cal I}_0^n\subset {\cal O}_Q$.
Let $Q^{[D_0.Y_0]}\to Q$
denote the dilatation 
for $Y_0\to Q$
and $D_0$.
We also define a dilatation
$Q^{[nD_0.Y_0^{(n)}]}\to Q$
for $Y_0^{(n)}\to Q$
and $nD_0$.

Since $Y_0$ is \'etale over $X$,
the scheme $Q^{[D_0.Y_0]}$
is smooth over $X$ 
by Example \ref{exQD}.1 and
equals its normalization
 $Q^{(D_0.Y_0)}$.
The canonical morphism
$Q^{[D_0.Y_0]}\to Q^{[nD_0.
Y_0^{(n)}]}$ is finite and
induces an isomorphism
$Q^{(D_0.Y_0)}\to Q^{(nD_0.
Y_0^{(n)})}$on the normalizations.

By the assumptions
${\cal J}_0^n=0$
and $nD_0\leqq D$,
we have 
${\cal I}_0^n
\subset 
{\cal I}+{\cal I}_D
\subset 
{\cal I}+{\cal I}_{nD_0}$.
Hence we have
a morphism
$Q^{[nD_0]}\to Q^{[nD_0.
Y_0^{(n)}]}$.
Further by $nD_0\leqq D$,
we obtain a morphism
$Q^{(D)}\to Q^{(D_0.Y_0)}$
of normalizations.

The dilatation $P^{[D]}$ of $P$
for the section $X\to P$ and $D$
is smooth over $X$ 
by Example \ref{exQD}.1 and
hence is equal to the normalization
$P^{(D)}$.
Since $Y^{(D)}\to X$
is \'etale
and since the diagram
$$
\begin{CD}
Y^{(D)}@>>>Q^{(D)}\\
@VVV@VVV\\
Y@>>>Q^{[D]}@>>>Q\\
@VVV@VVV@VVV\\
X@>>>P^{[D]}@>>> P
\end{CD}$$ 
is cartesian by Lemma \ref{lmFunY}.2,
the quasi-finite morphism
$Q^{(D)}\to P^{(D)}$
of normal schemes
is \'etale on a neighborhood 
$W\subset Q^{(D)}$ of $Y^{(D)}$
by \cite[Th\'eor\`eme (18.10.16)]{EGA4}.

The commutative diagram
$$\xymatrix{
Q^{(D)}\ar[r]\ar[d]&
Q^{(D_0.Y_0)}\ar[r]&
Q\ar[d]\\
P^{(D)}\ar[rr]&&
P}
$$
of schemes
defines a commutative diagram
$$\xymatrix{
{\cal O}_W\otimes
\Omega^1_{Q^{(D)}/X}&
{\cal O}_W\otimes
\Omega^1_{Q^{(D_0.Y_0)}/X}\ar[l]&
{\cal O}_W\otimes
\Omega^1_{Q/X}\ar[l]\\
{\cal O}_W\otimes
\Omega^1_{P^{(D)}/X}\ar[u]&&
{\cal O}_W\otimes
\Omega^1_{P/X}\ar[u]\ar[ll]}
$$
of locally free
${\cal O}_W$-modules.
Since $Q^{(D)}\to P^{(D)}$ is
\'etale on $W$,
the left vertical arrow
is an isomorphism.

Since $X\to X$ and $Y_0\to X$
are \'etale,
the lower horizontal arrow
(resp.\
the upper right horizontal arrow)
induces an isomorphism
${\cal O}_W\otimes
\Omega^1_{P/X}
\to
{\cal I}_D\cdot
{\cal O}_W\otimes
\Omega^1_{P^{(D)}/X}$
(resp.\ 
${\cal O}_W\otimes
\Omega^1_{Q/X}
\to
{\cal I}_{D_0}\cdot
{\cal O}_W\otimes\Omega^1_{Q^{(D_0.Y_0)}/X}$).
Hence ${\cal I}_{D-D_0}\cdot
{\cal O}_W\otimes\Omega^1_{Q/X}
=
{\cal I}_D\cdot
{\cal O}_W\otimes
\Omega^1_{Q^{(D_0.Y_0)}/X}$
is contained in the image
of 
${\cal O}_W\otimes\Omega^1_{P/X}$.
Or equivalently,
${\cal O}_W\otimes_{{\cal O}_Q}
\Omega^1_{Q/P}$
is annihilated by
${\cal I}_{D-D_0}$.
Hence its pull-back 
${\cal O}_{Y^{(D)}}\otimes
_{{\cal O}_Y}
\Omega^1_{Y/X}$
is also annihilated by
${\cal I}_{D-D_0}$.
\qed }

\section{Ramification}

\subsection{Ramification of quasi-finite schemes}\label{ssram}

Let $X$ be a normal noetherian
scheme and 
$D$ be an effective Cartier divisor of $X$.
Let $Y$ be a quasi-finite scheme 
over $X$
such that $D_Y=D\times_XY
\subset Y$ is a Cartier divisor.

Locally on $X$, there exists a smooth scheme
$Q$ over $X$
and a closed immersion $Y\to Q$ over $X$.
Then, by Proposition \ref{prFunQ}
and Corollary \ref{corFunQ}, 
the scheme $Y^{(D)}$ over $Y$
defined as $Y\times_{Q^{[D]}}Q^{(D)}$
is canonically independent of $Q$.
Hence a finite scheme $Y^{(D)}$ over $Y$ is
defined by patching.
Similarly, 
for a geometric point $\bar x$ above $x\in D$,
the set $\pi_0(Q^{(D)}_{\ \bar x})$
of connected components of
the geometric fiber 
is canonically independent of $Q$.

\begin{df}\label{dfF}
Let $X$ be a normal noetherian
scheme and 
$D$ be an effective Cartier divisor of $X$.
Let $Y$ be a quasi-finite scheme over $X$ such that $D_Y=D\times_XY
\subset Y$ is a Cartier divisor
and let $\bar Y$ be the normalization
of $Y$.
Let $\bar x$ be a geometric point above a point $x\in D$.

By taking a closed immersion
$Y\to Q$ to a smooth scheme $Q$ over $X$
defined on a neighborhood of $x$,
we define finite sets
$F_{\bar x}^D(Y/X)$
and $F_{\bar x}^{D+}(Y/X)$
by 
\begin{equation}
F_{\bar x}^D(Y/X)
=\pi_0(Q^{(D)}_{\ \bar x}),
\qquad
F_{\bar x}^{D+}(Y/X)
=Y^{(D)}_{\ \bar x}
\label{eqFD}
\end{equation}
equipped with canonical mappings
\begin{equation}
\xymatrix{
\bar Y_{\bar x}
\ar[r]^{\hspace{-5mm}\varphi^{D+}_{\bar x}}
\ar[d]_{\varphi^{D}_{\bar x}}
&
F_{\bar x}^{D+}(Y/X)
\ar[d]\ar[dl]
\\
F_{\bar x}^{D}(Y/X)
\ar[r]
&Y_{\bar x}
}
\label{eqFDc}
\end{equation}
induced by the morphisms
$$
\xymatrix{
\bar Y\ar[r]\ar[d]
&
Y^{(D)}\ar[d]\ar[dl]
\\
Q^{(D)}\ar[r]
&Q
}$$
in {\rm (\ref{eqYQD})}.
\end{df}

We consider a commutative diagram
\begin{equation}
\begin{CD}
Y'@>>> X'&\, \supset D'\supset\, & D\times_XX'
@<<<\bar x'\\
@VVV@VVV@VVV@VVV\\
Y@>>> X&\supset & D
@<<<\bar x
\end{CD}
\label{eqFunDY}
\end{equation}
of noetherian schemes.
We assume that
$X'$ is normal,
$D'\subset X'$ is an effective
Cartier divisor,
$Y'$ is quasi-finite
over $X'$ and that
$D'_{Y'}\subset Y'$ is
an effective Cartier divisor.
Then, the commutative diagram (\ref{eqFun[]})
induces a commutative diagram
\begin{equation}
\begin{CD}
\bar Y'_{\bar x'}
@>{\varphi^{D'+}_{\bar x'}}>>
F_{\bar x'}^{D'+}(Y'/X')
@>>>
F_{\bar x'}^{D'}(Y'/X')
@>>>
Y'_{\bar x'}\\
@VVV@VVV@VVV@VVV\\
\bar Y_{\bar x}
@>{\varphi^{D+}_{\bar x}}>>
F_{\bar x}^{D+}(Y/X)
@>>>
F_{\bar x}^{D}(Y/X)
@>>>
Y_{\bar x}.
\end{CD}
\label{eqFunDD'}
\end{equation}

For effective Cartier divisors
$D$ and $D'$ of a scheme $X$
defined by the ideal sheaves
${\cal I}_D,{\cal I}_{D'}\subset
{\cal O}_X$
and for $x\in D$,
we write $D<D'$ at $x$
if we have a strict inclusion
${\cal I}_{D,x}\supsetneqq
{\cal I}_{D',x}$.
If $X=X', Y=Y', \bar x=\bar x'$
and if $D<D'$ at the image $x$ of $\bar x$ as Cartier divisors, 
further we have an arrow
$F_{\bar x}^{D'}(Y/X)
\to
F_{\bar x}^{D+}(Y/X)$
making the two triangles
obtained by dividing the middle
square commutative
by Lemma \ref{lmFunD}.

\begin{pr}\label{prcoca}
Assume that
$Y\to X$ is quasi-finite, flat and
locally of complete intersection
and that the normalization
$\bar Y$ of $Y$
is \'etale over $X$.

{\rm 1.}
The arrows in
\begin{equation*}
\xymatrix{
\bar Y_{\bar x}
\ar[r]^{\hspace{-5mm}\varphi^{D+}_{\bar x}}
\ar[d]_{\varphi^{D}_{\bar x}}
&
F_{\bar x}^{D+}(Y/X)
\ar[d]\ar[dl]
\\
F_{\bar x}^{D}(Y/X)
\ar[r]
&Y_{\bar x}
}
\leqno{(\ref{eqFDc})}
\end{equation*}
are surjections.

{\rm 2.}
Let $Y'\to Y$
be a surjective morphism 
locally of complete intersection
of quasi-finite and flat schemes
over $X$.
Assume that the normalization
$\bar Y'$ of $Y'$
is \'etale over $X$.
Then, the diagram
\begin{equation}
\begin{CD}
\bar Y'_{\bar x}
@>{\varphi^{D+}_{\bar x}}>>
F^{D+}_{\bar x}(Y'/X)
@>>>
F^{D}_{\bar x}(Y'/X)
@>>>
Y'_{\bar x}
\\
@VVV
@VVV@VVV@VVV\\
\bar Y_{\bar x}
@>{\varphi^{D+}_{\bar x}}>>
F^{D+}_{\bar x}(Y/X)
@>>>
F^{D}_{\bar x}(Y/X)
@>>>
Y_{\bar x}
\end{CD}
\label{eqcoca}
\end{equation}
is a cocartesian diagram of surjections.
\end{pr}

\proof{
By replacing $X$ by the strict localization $X_{(\bar x)}$,
we may assume that
$\bar x\to X$ is a closed immersion
and that $Y\to X$ is finite.

1.
By Lemma \ref{lmPQ}.2,
we may assume that
there exists smooth
schemes $P$ and $Q$ over $X$ 
and a cartesian diagram
$$\begin{CD}
Y@>>> Q\\
@VVV@VVV\\
X@>>>P
\end{CD}$$
of schemes over $X$
such that the horizontal arrows
are closed immersions
and that the vertical arrows
are quasi-finite and flat.
We verify that the diagram
\begin{equation}
\begin{CD}
Q^{(D)}_{\ \bar x}@>>>Q^{(D)}@<<< Y^{(D)}
@<<< \bar Y\\
@VVV
\hspace{-10mm}
\square
\hspace{7mm}
@VVV
@VVV@VVV\\
P^{(D)}_{\ \bar x}@>>>P^{(D)}@<<< X
@= X
\end{CD}
\label{eqPDX}
\end{equation}
satisfies the assumptions
in Corollary \ref{corCA}.
Since $P^{[D]}\to X$ is smooth,
we have $P^{(D)}=P^{[D]}$.
By Lemma \ref{lmFunY}.2,
the diagram
$$\begin{CD}
Q@<<<Q^{[D]}@<<< Y\\
@VVV@VVV@VVV\\
P@<<<P^{[D]}@<<< X
\end{CD}$$
is cartesian.
Hence the middle square in
(\ref{eqPDX}) is also cartesian.

The diagram (\ref{eqPDX}) satisfies
the finiteness
assumption in Corollary \ref{corCA},
by Lemma \ref{lmFunY}.1.
Since $X=X_{(\bar x)}$ is strictly local,
the assumption that
the canonical mapping
${\bar x}\to 
\pi_0(\bar X)$ is a bijection is satisfied.
Since $P^{(D)}_{\ \bar x}$ 
is a vector space over $\bar x$
and is connected,
the mapping
${\bar x}\to 
P^{(D)}_{\ \bar x}\cap X
\to \pi_0(P^{(D)}_{\ \bar x})$ 
are bijections of sets consisting
of single elements.
We may assume that the finite \'etale morphism
$\bar Y\to X$
is surjective since if otherwise the assertion is trivial.
Hence by
Corollary \ref{corCA}.2
(resp.\ 3),
the mapping $\bar Y_{\bar x}\to
Y^{(D)}_{\ \bar x}=
F^{D+}_{\bar x}(Y/X)$
(resp.\
$F^{D+}_{\bar x}(Y/X)
=Y^{(D)}_{\ \bar x}
\to 
\pi_0(Q^{(D)}_{\ \bar x})
=
F^{D}_{\bar x}(Y/X)$)
is surjective.

Similarly, applying 
Corollary \ref{corCA}.2
to the diagram
$$\begin{CD}
Y_{\bar x}@>>>Y@<<< Y
@<<< \bar Y\\
@VVV
\hspace{-10mm}
\square
\hspace{7mm}
@VVV
\hspace{-10mm}
\square
\hspace{7mm}
@VVV@VVV\\
\bar x@>>>X@<<< X
@= X,
\end{CD}$$
we see that
$\bar Y_{\bar x}\to
Y_{\bar x}$ is a surjection.

2.
By Lemma \ref{lmPQ}.2,
we may assume that
there exists smooth
schemes $Q$ and $Q'$ over $X$ 
and a cartesian diagram
$$\begin{CD}
Y'@>>> Q'\\
@VVV@VVV\\
Y@>>>Q
\end{CD}$$
of schemes over $X$
such that the horizontal arrows
are closed immersions
and that the vertical arrows
are quasi-finite and flat.

We verify that the diagram
$$\begin{CD}
Q^{\prime (D)}_{\ \bar x}
@>>>Q^{\prime (D)}@<<< Y^{\prime (D)}
@<<< \bar Y'\\
@VVV
\hspace{-10mm}
\square
\hspace{7mm}
@VVV
@VVV@VVV\\
Q^{(D)}_{\ \bar x}@>>>Q^{(D)}@<<< Y^{(D)}
@<<< \bar Y
\end{CD}$$
satisfies the assumptions in
Corollary \ref{corCA}.
The middle square is cartesian by
Lemma \ref{lmFunY}.2.
The finiteness
assumption in Corollary \ref{corCA}
is satisfied
by Lemma \ref{lmFunY}.1.
Since the finite \'etale
covering $\bar Y\to X$ is split
and $X$ is connected,
the assumption that
the canonical mapping
$\bar Y_{\bar x}\to 
\pi_0(\bar Y)$ is a bijection is satisfied.
By 1,
$\bar Y_{\bar x}
\to Y^{(D)}_{\ \bar x}
\to \pi_0(Q^{(D)}_{\ \bar x})$
are surjective.
We may assume that $Y$ and $Y'$ are finite over 
$X$. Since $Y'\to Y$ is surjective,
the morphism $\bar Y'\to \bar Y$
of finite \'etale schemes over $X$
is also surjective.
Hence by
Corollary \ref{corCA}.2
(resp.\ 3),
the right square
(resp.\ the middle square)
of (\ref{eqcoca})
is a cocartesian diagram of surjections.

Similarly, applying 
Corollary \ref{corCA}.2
to the diagram
$$\begin{CD}
Y'_{\bar x}@>>>
Y'@<<< Y'
@<<< \bar Y'\\
@VVV
\hspace{-10mm}
\square
\hspace{7mm}
@VVV
\hspace{-10mm}
\square
\hspace{7mm}
@VVV@VVV\\
Y_{\bar x}@>>>Y@<<< Y
@<<< \bar Y,
\end{CD}$$
we see that
the big rectangle in (\ref{eqcoca})
is a cocartesian diagram of surjections.
\qed }

\begin{cor}\label{coretaYD}
Assume that $Y\to X$ is locally of
complete intersection.
Let $P$ and $Q$ be smooth
schemes over $X$ and let 
$$\begin{CD}
Y@>>> Q\\
@VVV
\hspace{-10mm}
\square
\hspace{7mm}
@VVV\\
X@>>>P
\end{CD}$$
be a cartesian diagram
of schemes over $X$
such that the horizontal arrows
are closed immersions
and that the vertical arrows
are quasi-finite and flat.
Then, 
the mapping 
$\bar Y_{\bar x}\to 
F^{D+}_{\bar x}(Y/X)$
is an injection on the inverse image
of $y\in Y$ if and only if
$Q^{(D)}\to P^{(D)}$
is \'etale on the inverse image of
$y$ by $Y^{(D)}\to Y$.
\end{cor}

\proof{
Since the assertion is \'etale local,
we may assume that $Y\to X$ 
and $Q\to P$ are finite
and that $y$ is the unique point
of the inverse image of $x$.
Then, by Proposition \ref{prcoca}.1,
$\bar Y_{\bar x}\to Y^{(D)}_{\bar x}
=F^{D+}_{\bar x}(Y/X)
\subset Q^{(D)}_{\bar x}$
is a bijection of finite sets.
Hence 
$Q\to P$ is \'etale at $x$
by \cite[Th\'eor\`eme (18.10.16)]{EGA4}.
\qed }

\begin{df}\label{dfram}
Let $X$ be a normal noetherian
scheme and
$U\subset X$ be a dense open subscheme.
Let $Y$ be a quasi-finite scheme over $X$
such that $V=U\times_XY\to U$
is \'etale.
Let $D$ be an effective Cartier divisor of
$X$ such that $U\cap D$ is empty
and that $D_Y=D\times_XY$
is an effective Cartier divisor.

{\rm 1.}
For $x\in D$,
we consider the following condition
on $X,Y$ and $D$:

{\rm (RF)}  
There exist an open neighborhood
$W$ of $x\in X$,
a smooth scheme $Q$ 
over $W$ and a closed immersion
$Y\times_XW\to Q$ 
such that the normalization
$\bar Y$ of $Y$ is \'etale over $W$
and that the normalization
$Q^{(D)}$ of the dilatation
$Q^{[D]}$ is flat and reduced
over $W$.

If the condition {\rm (RF)} 
is satisfied at every $x\in D$,
we say that
$Y$ over $X$ satisfies the condition
{\rm (RF)} for $D$.

{\rm 2}.
Let $x\in D$ and assume that
$Y$ over $X$ satisfies the condition
{\rm (RF)} for $D$ at $x$.

Let $y$ be a point of 
$\bar Y\times_Xx
\subset
\bar Y\times_XD$.
We say that the ramification of
$Y\to X$ is bounded by $D$ 
(resp.\ by $D+$) at $y$,
if the mapping
$\varphi^{D}_{\bar x}\colon \bar Y_{\bar x}
\to F_{\bar x}^D(Y/X)$
(resp.\ 
$\varphi^{D+}_{\bar x}
\colon \bar Y_{\bar x}
\to F_{\bar x}^{D+}(Y/X)$)
is an injection on the inverse
image of $y$.

We say that the ramification of
$Y\to X$ is bounded by $D$ 
(resp.\ by $D+$) at $x$,
if the mapping
$\varphi^{D}_{\bar x}\colon \bar Y_{\bar x}
\to F_{\bar x}^D(Y/X)$
(resp.\ 
$\varphi^{D+}_{\bar x}
\colon \bar Y_{\bar x}
\to F_{\bar x}^{D+}(Y/X)$)
is an injection.
\end{df}

If ramification is bounded by $D$,
it is bounded by $D+$.
We show that the condition
(RF) is independent of the
choice of $Q$.

\begin{lm}\label{lmRF}
Let $X$ be a normal noetherian
scheme and
$U\subset X$ be a dense open subscheme.
Let $Y$ be a quasi-finite scheme over $X$
such that $V=U\times_XY\to U$
is \'etale.
Let $D$ be an effective Cartier divisor of
$X$ such that $U\cap D$ is empty
and that $D_Y\subset Y$ is an
 effective Cartier divisor.
Let $x\in D$.

{\rm 1.}
Assume that $Y$ over $X$
satisfies {\rm (RF)} for $D$
at $x$.
Let $W\subset X$ be an
open neighborhood of $x$,
$Q$ be a smooth scheme over $W$
and $Y\times_XW\to Q$
be a closed immersion.
Then, there exists an
open neighborhood $W'\subset W$
of $x$, such that
$(Q\times_WW')^{(D\times_XW')}
\to W'$ is flat and reduced.

{\rm 2.}
Let $X'\to X$ be a morphism
of normal noetherian scheme
such that $U'=U\times_XX'$
is a dense open subscheme
and that
$D'_{Y'}=
D_Y\times_XX'
\subset Y'=Y\times_XX'$
is an effective Cartier divisor.
Let $x'$ be a point of
$D'=D\times_XX'$ above $x$.
We consider the following conditions:

{\rm (1)}
$Y$ over $X$
satisfies {\rm (RF)} for $D$
at $x$.

{\rm (2)}
$Y'$ over $X'$
satisfies {\rm (RF)} for $D'$
at $x'$.

\noindent
We have {\rm (1)}$\Rightarrow${\rm (2)}.
Conversely,
if $X'\to X$
is smooth at $x'$, 
we have {\rm (2)}$\Rightarrow${\rm (1)}.
\end{lm}

\proof{
1. Set $D_W=D\times_XW$.
After shrinking $W$ if necessary,
we may assume that
there exist a smooth scheme
$Q_0$ over $W$
and a closed immersion
$Y\times_XW\to Q_0$
such that
$Q_0^{(D_X)}
\to W$ is flat and reduced.
Since 
$Q^{(D_W)}
\gets
(Q\times_WQ_0)^{(D_W)}
\to
Q_0^{(D_W)}$
are smooth by
Proposition \ref{prFunQ},
the assertion follows.

2.
{\rm (1)}$\Rightarrow${\rm (2)}:
This follows from Lemma \ref{lmFunX}.

{\rm (2)}$\Rightarrow${\rm (1)}:
After shrinking
$X'$ if necessary,
we may assume that $X'\to X$ is smooth.
Let $W$ be an open neighborhood
of $x$, let
$Y\times_XW\to Q$ be 
a closed immersion to a smooth
scheme $Q$ over $W$
and $W'=W\times_XX'$.
Then the morphism
$(Q\times_WW')^{(D'\times_{X'}W')}
\to
Q^{(D\times_XW)}\times_WW'$
is an isomorphism 
by Lemma \ref{lmFunX}.
Hence the assertion follows.
\qed }

\begin{lm}\label{lmbdy}
Let $X$ be a normal noetherian
scheme and
$U\subset X$ be a dense open subscheme.
Let $Y$ be a quasi-finite scheme over $X$
such that $V=U\times_XY\to U$
is \'etale.
Let $D\subset D'$ be effective Cartier divisors of
$X$ such that $U\cap D'$ is empty
and that $D'_Y\subset Y$
is an effective Cartier divisor.
Let $x\in D$ 
and assume that $Y$ over $X$
satisfies {\rm (RF)} for $D$
and $D'$ at $x$.

Let $y\in Y$ be a point above $x$.
If the ramification of $Y$ over $X$ is
bounded by $D+$ at $y$
and if $D<D'$ at $x$,
then the ramification of $Y$ over $X$ 
is bounded by $D'$ at $y$.
\end{lm}

\proof{
It follows from Lemma \ref{lmFunD}.
\qed }

\begin{lm}\label{lmramD}
Let $X$ be a normal noetherian
scheme and
$U\subset X$ be a dense open subscheme.
Let $Y$ be a quasi-finite scheme over $X$
such that $V=U\times_XY\to U$
is \'etale.
Let $D$ be an effective Cartier divisor of
$X$ such that $U\cap D$ is empty
and that $D_Y\subset Y$
is an effective Cartier divisor.
Assume that $Y$ over $X$
satisfies the condition {\rm (RF)}
for $D$.

Let $S\subset D_Y$ 
(resp.\ $S^+\subset D_Y$)
denote the subset consisting of points $y\in D_Y$
where the ramification of
$Y\to X$ is bounded by $D$
(resp.\ by $D+$).

{\rm 1.}
We have $S\subset S^+$.

{\rm 2.}
The subset $S\subset D_Y$ 
is closed
and 
the subset $S^+\subset D_Y$ 
is open.
\end{lm}

\proof{1.
It follows from the commutative
diagram (\ref{eqFD}).

2.
By Lemma \ref{lmA}
applied to $\bar Y\to Q^{(D)}$,
we see that $S$ is closed.
Similarly, by Lemma \ref{lmB}
applied to $\bar Y\to Y^{(D)}$
we see that $S^+$ is open.
\qed }

\begin{pr}\label{prstD}
Let $X$ be a normal noetherian
scheme and
$U\subset X$ be a dense open subscheme.
Let $Y$ be a quasi-finite scheme over $X$
such that $V=U\times_XY\to U$
is \'etale.
Let $D$ be an effective Cartier divisor of
$X$ such that $U\cap D$ is empty
and that $D_Y\subset Y$
is an effective Cartier divisor.

Let $C$ be a semi-stable curve
over $X$ such that
$C_U=C\times_XU\to U$ is smooth.
Let $x\in X$ be a point of $D$
and $z\in C$ be a singular
point of the fiber $C_x$.
Assume that there exist
two irreducible components
$C_1$ and $C_2$ of the fiber $C_x$
meeting at $z$ and 
let $\zeta_1$ and $\zeta_2$
be their generic points.
Let $D_1\subset D_2$ 
be effective Cartier divisors on $X$
and let $\tilde D\subset C$ 
be an effective Cartier divisor 
such that 
$D_1<D_2$ at $x$ and that 
$\tilde D=D_i\times_XC
=D_{i,C}$
on a neighborhood of $\zeta_i$
for $i=1,2$.

Assume that $Y_C=Y\times_XC$ 
over $C$
satisfies the condition {\rm (RF)}
for $\tilde D$ at $z$.

{\rm 1.}
$Y$ over $X$
satisfies the condition {\rm (RF)}
for $D_1$ and $D_2$
at $x$.

{\rm 2.}
We have a commutative diagram
\begin{equation}
\xymatrix{
F_{\bar x}^{D_2+}(Y/X)
\ar[d]\ar[r]
&
F_{\bar z}^{\tilde D+}(Y_C/C)
\ar[d]\ar[r]
&
F_{\bar x}^{D_1+}(Y/X)
\ar[d]\\
F_{\bar x}^{D_2}(Y/X)
\ar[ru]\ar[r]
&
F_{\bar z}^{\tilde D}(Y_C/C)
\ar[ru]\ar[r]
&
F_{\bar x}^{D_1}(Y/X).
}
\label{eqsst}
\end{equation}

{\rm 3.}
The lower left horizontal arrow
$F_{\bar x}^{D_2}(Y/X)\to
F_{\bar z}^{\tilde D}(Y_C/C)$
in {\rm (\ref{eqsst})}
is an injection.
The upper right horizontal arrow
$F_{\bar z}^{\tilde D+}(Y_C/C)
\to
F_{\bar x}^{D_1+}(Y/X)$
in {\rm (\ref{eqsst})}
is an injection on the image
of $\bar Y_{\bar x}$.
\end{pr}

\proof{
1.
Since $\zeta_1$ and $\zeta_2$
are contained in any open neighborhood of $z$,
the scheme $Y_C$ over $C$
satisfies  {\rm (RF)}
for $\tilde D$
at $\zeta_1$ and $\zeta_2$.
Since $C\to X$ is smooth at
$\zeta_1$ and $\zeta_2$,
the scheme $Y$ over $X$
satisfies  {\rm (RF)}
for $D_1$ and $D_2$
at $x$ by Lemma \ref{lmRF}.2.

2.
Let $D_{1,C}$ and $D_{2,C}$
be the pull-backs of
$D_1$ and $D_2$ to $C$.
Then, we have
$D_{1,C}<\tilde D<D_{2,C}$
at $z$. Hence by (\ref{eqFunDD'})
with the slant arrow added,
we obtain a commutative diagram
\begin{equation}
\xymatrix{
F_{\bar z}^{D_{2,C}+}(Y_C/C)
\ar[d]\ar[r]
&
F_{\bar z}^{\tilde D+}(Y_C/C)
\ar[d]\ar[r]
&
F_{\bar z}^{D_{1,C}+}(Y_C/C)
\ar[d]\\
F_{\bar z}^{D_{2,C}}(Y_C/C)
\ar[ru]\ar[r]
&
F_{\bar z}^{\tilde D}(Y_C/C)
\ar[ru]\ar[r]
&
F_{\bar z}^{D_{1,C}}(Y_C/C).
}
\label{eqsstC}
\end{equation}
Since $Y$ over $X$
satisfies  {\rm (RF)}
for $D_1$ and $D_2$
at $x$ by 1,
the pull-back defines
canonical isomorphisms
from the left and right
columns of (\ref{eqsst}) to
those of (\ref{eqsstC})
by Lemma \ref{lmFunX}.
Thus we obtain (\ref{eqsst}).

3.
By functoriality of cospecialization
mappings, we obtain a commutative
diagram
\begin{equation}
\xymatrix{
F_{\bar \zeta_2}^{D_{2,C}}(Y_C/C)
\ar[d]&
F_{\bar z}^{D_{2,C}}(Y_C/C)
\ar[l]_{\rm cosp.}\ar[d]&
F_{\bar x}^{D_2}(Y/X)
\ar[l]\ar[ld]\\
F_{\bar \zeta_2}^{\tilde D}(Y_C/C)&
F_{\bar z}^{\tilde D}(Y_C/C).
\ar[l]_{\rm cosp.}&
}
\label{eqsstc}
\end{equation}
By Lemma \ref{lmFunX} and by $\tilde D=D_{2,C}$
at $\zeta_2$,
the composition
$F_{\bar x}^{D_2}(Y/X)\to
F_{\bar \zeta_2}^{\tilde D}(Y_C/C)$
is a bijection.
Hence $F_{\bar x}^{D_2}(Y/X)\to
F_{\bar z}^{\tilde D}(Y_C/C)$
is injective.

Since the second assertion is \'etale local
on $X$, we may assume that
$Y\to X$ is finite.
By functoriality of specialization
mappings, we obtain a commutative
diagram
$$
\xymatrix{
&
F_{\bar z}^{\tilde D+}(Y_C/C)
\ar[d]\ar[ld]&
F_{\bar \zeta_1}^{\tilde D+}(Y_C/C)
\ar[l]_{\rm sp.}\ar[d]&
\bar Y_{\bar x}\ar[dl]\ar[l]
\\
F_{\bar x}^{D_1+}(Y/X)
\ar[r]&
F_{\bar z}^{D_{1,C}+}(Y_C/C)&
F_{\bar \zeta_1}^{D_{1,C}+}(Y_C/C).
\ar[l]_{\rm sp.}&&
}
$$
Since the composition
$F_{\bar \zeta_1}^{\tilde D+}(Y_C/C)
\to
F_{\bar z}^{D_{1,C}+}(Y_C/C)$
is a bijection,
the vertical arrow
$F_{\bar z}^{\tilde D+}(Y_C/C)
\to
F_{\bar z}^{D_{1,C}+}(Y_C/C)$
is an injection on the image
of $\bar Y_{\bar x}$.
Hence the assertion follows.
\qed }

\subsection{Ramification and valuations}

For a valuation ring $A\subset K$,
let $v\colon K^\times\to
 \Gamma=K^\times/A^\times$
denote the valuation.

\begin{df}\label{dfvale}
Let $X$ be a normal separated
noetherian scheme,
$U\subset X$ be 
a dense open subscheme
and $A$ be a valuation ring.
We say that a morphism
$T={\rm Spec}\, A\to X$ 
is $U$-{\em external} if $T\times_XU$
consists of a single point $t$.

For a morphism
$T={\rm Spec}\, A\to X$
and an effective Cartier divisor $D\subset X$,
let $v(D)\in \Gamma$
denote the valuation
$v(f)$ of a non-zero divisor $f$
defining $D\subset X$ 
on a neighborhood of
the image of $T$.
\end{df}

Let $\widetilde X=
\varprojlim X'$
be the inverse limit of proper
schemes $X'\to X$
such that
$U'=U\times_XX'\to U$
is an isomorphism.
Then, 
points of $\widetilde X
\sm U$ correspond 
bijectively to the inverse limits of the images 
of the closed points
by the liftings of $U$-external morphisms
$T\to X$ defined by
valuation rings of the residue
fields of points of $U$
by \cite[5.4]{FK}.

\begin{lm}\label{lmcof}
Let $X$ be a normal noetherian scheme,
$U\subset X$ be a dense open subscheme,
$t\in U$ be a point,
$A\subsetneqq K=k(t)$ be
a valuation ring
and $T={\rm Spec}\, A
\to X$ be a $U$-external morphism.

{\rm 1.}
Let $g\in \Gamma(U',{\cal O}^\times_{U'})$
be an invertible function
defined on an open neighborhood 
$U'\subset U$ of
$t\in U$
such that $v(g)=\gamma\geqq 0$.
Then, there exist a normal scheme 
$X'$ of finite type over $X$
such that $U\times_XX'=U'$,
that $g$ is extended
to a non-zero divisor on $X'$
defining an effective Cartier divisor
$R'\subset X'$
and that $U'=X'\sm D'$ is the complement
of an effective Cartier divisor $D'
\subset X'$ and a $U'$-external morphism $T\to X'$ lifting
$T\to X$ and $v(R')=\gamma$.

{\rm 2.}
Let $K'$ be a finite separable
extension of $K=k(t)$ and
$A'\subsetneqq K'$
be a valuation ring such
that $A'\cap K=A$.
Set $T'={\rm Spec}\, A'$
and let $\gamma>0$
be a positive element of the value
group $\Gamma'$ of $A'$.
Then, there exist a commutative diagram
$$\begin{CD}
U'@>>> X'@<<<T'\\
@VVV
\hspace{-10mm}
\square
\hspace{7mm}
@VVV@VVV\\
U@>>> X@<<<T
\end{CD}$$
of schemes,
a point $t'\in U'$ above $t$,
an isomorphism $K'\to k(t')$ over $K$
and an effective Cartier divisor $R'$
of $X'$
satisfying the following conditions {\rm (i)-(iv)}:

{\rm (i)} $X'$ is a normal scheme
of finite type over $X$,

{\rm (ii)} The left square is cartesian and
$U'$ is a dense open subscheme
of $X'$ \'etale over $U$.

{\rm (iii)} $T'\to X'$ is a $U'$-external morphism
extending $t'\to U'$.

{\rm (iv)} $R'\cap U'=\varnothing$
and $v'(R')=\gamma$.

{\rm 3.}
Let $$
\xymatrix{
&U'\ar[r]\ar[ddl]&
X'\ar[ddl]&
T'\ar[ddl]\ar[l]&
\bar x'\ar[ddl]\ar[l]
\\
U_1\ar[r]\ar[d]&
X_1\ar[d]&
T_1\ar[d]\ar[l]&
\bar x_1\ar[d]\ar[l]&
\\
U\ar[r]&
X&
T\ar[l]&
\bar x\ar[l]&
}$$
be a commutative diagram,
$t_1\in U_1$ and $t'\in U'$ be
points above $t\in U$
and $R_1\subset X_1$
and $R'\subset X'$ be
effective Cartier divisors
satisfying the following conditions
{\rm (i)-(iv)}:

{\rm (i)} $X_1$ and $X'$
are normal noetherian schemes
and $X_1\to X$ is of finite type.

{\rm (ii)}
The left square and the
left parallelogram are cartesian
and $U_1\to U$ is \'etale.
The open subschemes
$U_1\subset X_1$
and $U'\subset X'$ are dense.

{\rm (iii)}
$T_1={\rm Spec}\, A_1$
and $T'={\rm Spec}\, A'$
for valuation rings 
$A_1\subsetneqq K_1=k(t_1)$
and
$A'\subsetneqq K'=k(t')$
satisfying
$A_1\cap K=A'\cap K=A$.
The morphism
$T_1\to X_1$ is $U_1$-external
and $T'\to X'$ is $U'$-external.

{\rm (iv)}
$R_1\cap U_1$
and $R'\cap U'$ are empty
and we have
$v_1(R_1)\leqq
v'(R')$ in $\Gamma'_{\mathbf Q}$.

\noindent
Then, there exist a commutative
diagram
$$\begin{CD}
U'_1@>>> X'_1@<<<T'_1@<<<\bar x'_1\\
@|@VVV@VVV@VVV\\
U'\times_UU_1@>>> 
X'\times_XX_1@<<<
T'\times_TT_1@<<<
\bar x'\times_{\bar x}\bar x_1
\end{CD}$$ 
and $t'_1\in U'_1$ 
above $t$
satisfying the following conditions
{\rm (i)-(iv)}:

{\rm (i)}
$X'_1$ is a normal scheme of
finite type over $X'$.

{\rm (ii)}
$U'_1$ is $U'\times_{X'}X'_1$ and
is a dense open subscheme of $X'_1$.

{\rm (iii)}
$T'_1={\rm Spec}\, A'_1$
for a valuation ring $A'_1\subset K'_1=k(t'_1)$.

{\rm (iv)}
For the pull-backs
$R'_1=R_1\times_{X_1} X'_1$ and
$R'_2=R'\times_{X'} X'_1$,
we have
$R'_1\leqq R'_2$.
\end{lm}

\proof{
1.
Let $Z$ and $Z'$ be closed subschemes
such that
$U=X\sm Z$ and
$U'=X\sm Z'$.
By replacing $X$ by the normalization
of the blow-up at $Z$ and at $Z'$
and by the valuative criterion of properness,
we may assume that
$U=X\sm D$ and
$U'=X\sm D'$ are the complements
of effective Cartier divisors $D, D'\subset X$.

Let $x\in X$ be the image 
of the closed point of $T$
and let $W={\rm Spec}\, B\subset X$
be an open neighborhood of $x$
such that $D\cap W,D'\cap W$ are
principal divisors defined by
$f,f'\in B$.
Then we have $U'\cap W=
W\sm D'\cap W={\rm Spec}\, B[1/f']$.
Set $g=h/f^{\prime n}\in B[1/f']$.
The function $g$ and hence also
$h\in B$ are also invertible on $U'\cap W$.
Set $\alpha=v(f),
\alpha'=v(f')\in \Gamma$.
Since $T\to X$ is $U$-external,
we have $\Gamma^{+}[1/\alpha]
=\Gamma$.
Hence after replacing $f$
by its power,
we may assume that
$\alpha'\leqq \alpha$.

Let $W'\to W$ be the
normalization of the blow-up
at the ideals $(f^{\prime n}, h)$
and $(f,f')$.
Since $f,f'$ and $h$ are invertible
on $U'\cap W$,
the morphism
$W'\to W$ induces 
an isomorphism 
$U'\times_{X'}W'\to U'\cap W$.
Since $W'\to W$ is proper,
the morphism $T\to W$
is uniquely lifted to $T\to W'$.
Since the generic point
$t\in T$ is the unique
point of $U\times_XT
\supset (U'\times_{X'}W')\times_{W'}T$,
the morphism
$T\to W'$ is $U'$-external.

Let $x'\in W'$ be the image
of the closed point of $T$.
Since the ideals $(f^{\prime n}, h),
(f,f')\subset
{\cal O}_{W,x'}$
are principal ideals
and since
$v(h)\geqq v(f^{\prime n})$ and 
$v(f')\geqq v(f)$,
there exists an open 
neighborhood $X'$ of
$x'\in W$ such that
$U'\subset X'$ where
we have inclusions
$(f^{\prime n})\supset (h)$ and
$(f)\supset (f')$.
Then, $g=h/f^{\prime n}$ defines
a Cartier divisor $R'$ on $X'$
satisfying $R'\cap U'=\varnothing$
and $v(R')=\gamma$.
We also have an inclusion
$U\times_XX'
=X'\sm D\times_XX'
\subset
X'\sm D'\times_XX'=
U'\times_XX'=U'$.
Since the other inclusion is
obvious,
we have
$U'=U\times_XX'$.

2.
We may take an \'etale
scheme $U_1\to U$
such that $t'={\rm Spec}\, K'
=t\times_UU_1$
and a finite scheme
$X_1\to X$ containing
$U_1$ as a dense open scheme.
After shrinking $U_1$ 
if necessary,
we may take an invertible function
$g\in \Gamma(U_1,
{\cal O}_{U_1}^\times)$
such that $\gamma=v'(g)$.
Since $T'$ is a localization of
the normalization of
$T\times_XX_1$,
the morphism $t'\to U_1\subset X_1$
is uniquely extended to $T'\to X_1$.

Then, by 1 applied
to the open subschemes
$U_1\subset U\times_XX_1
\subset X_1$,
to the morphism
$T'\to X_1$
and to the invertible function
$g\in \Gamma(U_1,{\cal O}^\times_{U_1})$,
the assertion follows.

3.
Let $T_{(\bar x)},
T_{1,(\bar x_1)}$
and $T'_{(\bar x')}$
denote the strict localizations.
We take a point 
$\tilde t'_1\in 
T'_{(\bar x')}\times_{T_{(\bar x)}}
T_{1,(\bar x_1)}$
above the generic point of
$T'_{(\bar x')}$.
Then the normalization 
$\tilde T'_1$
of
$T'_{(\bar x')}$
in $\tilde t'_1$ is
$\tilde T'_1={\rm Spec}\, 
A^{\prime sh}_1$
for a strictly local valuation ring
$A^{\prime sh}_1$.
Let $t'_1\in 
t'\times_tt_1\subset
T'\times_TT_1$
be the image of $\tilde t'_1$
and set $K'_1=k(t'_1)$ and
$A'_1=A^{\prime sh}_1
\cap K'_1$.
Let $T'_1={\rm Spec}\, A'_1$
and $\bar x'_1$ be the geometric
point of $T'_1$
defined by a geometric
closed point of $\tilde T'_1$.

Let $X'_0$ be the normalization of
$X'\times_XX_1$ in
$U'_1=U'\times_UU_1$.
Define effective Cartier divisors
of $X'_0$ by 
$R'_{0,1}=R_1\times_{X_1}X'_0$
and
$R'_{0,2}=R'\times_{X'}X'_0$.
Let $\bar X'_1\to X'_0$ be the 
normalization of the blow-up
at $R'_0\cap R'_1
=R'_0\times_{X'_0}R'_1$
and
define effective Cartier divisors
of $\bar X'_1$ by 
$\bar R'_1=R_1\times_{X_1}\bar X'_1$
and
$\bar R'_2=R'\times_{X'}\bar X'_1$.
Since $\bar X'_1\to X'\times_XX_1$
is proper,
the morphism
$t'_1\to
t'\times_tt_1\subset
U'\times_UU_1$
is uniquely lifted to
$T'_1\to \bar X'_1$
by the valuative criterion
of properness.

Let $x'_1\in \bar X'_1$
be the image of 
the closed point of $T'_1$.
The intersection
$\bar R'_1\cap \bar R'_2\subset 
\bar X'_1$
is the exceptional divisor
and hence is an effective
Cartier divisor.
Since $v'_1(\bar R'_1)
\leqq v'_1(\bar R'_2)$,
on an open neighborhood
$X'_1\subset \bar X'_1$
of $x'_1$,
we have
$\bar R'_1\cap \bar R'_2
=\bar R'_1\leqq 
\bar R'_2$
by Nakayama's lemma.
\qed }

\medskip

Let $X$ be a normal noetherian scheme and
$U\subset X$ be a dense open 
subscheme. 
Let $t\in U$ and
$T={\rm Spec}\, A\to X$
be a $U$-external morphism
defined by a valuation ring
$A\subsetneqq K=k(t)$
of the residue field at a point $t\in U$.
Let $\bar x$ and $\bar t$
be geometric points of $T$
supported on the closed 
point and on the generic point
respectively.
Recall that $T_{(\bar x)}$
denotes the strict localization
and that
a specialization
$\bar x\gets \bar t$ is
a morphism $T_{(\bar x)}
\gets \bar t$ of schemes.

Let $A'$ be a valuation ring
and $T'={\rm Spec}\, A'\to T$ be
a faithfully flat morphism.
We identify $\Gamma$ as a subgroup
of the value group $\Gamma'$ of $A'$
by the canonical injection $\Gamma\to \Gamma'$.
Let $\bar x'$ and $\bar t'$ 
be geometric points
of $T'$ above
$\bar x$ and $\bar t$ respectively.
We say that a specialization
$\bar x'\gets \bar t'$ is a lifting of
$\bar x\gets \bar t$ if the diagram
$$\begin{CD}
\bar x'@>>> T'@<<< 
T'_{(\bar x')}@<<< \bar t'\\
@VVV@VVV@VVV@VVV\\
\bar x@>>> T@<<< 
T_{(\bar x)}@<<< \bar t
\end{CD}$$
is commutative.

We consider a commutative
diagram
\begin{equation}
\begin{CD}
X'@<<< T'\\
@VVV@VVV\\
X@<<< T
\end{CD}
\label{eqXT'}
\end{equation}
of schemes equipped with an
effective Cartier divisor
$R'\subset X'$ and a lifting
$\bar x'\gets \bar t'$ to $T'$
of the specialization $\bar x\gets \bar t$
satisfying the following conditions (i)-(iii):

(i) $X'$ is a normal noetherian scheme
of finite type over $X$
such that
$U'=U\times_XX'\subset X'$ is a 
dense open subscheme
\'etale over $U$.

(ii) $T'={\rm Spec}\, A'\to X'$
is a $U'$-external morphism
defined by a valuation ring
$A'\subsetneqq K'=k(t')$
of the residue field at a point $t'\in U'$
above $t$ such that
$A'\cap K=A$.

(iii)
$R'\cap U'=\varnothing$ and 
$v'(R')=\gamma$
in the value group $\Gamma'$
of $A'$.

For elements $\alpha\leqq \beta$
of a totally ordered group
$\Gamma$,
let $(\alpha,\beta)_\Gamma
\subset \Gamma$
denote the subset
$\{\gamma\in 
\Gamma\mid
\alpha<\gamma<\beta\}$.
Similarly,
we define
$(\alpha,\beta]_\Gamma,
(\alpha,\infty)_\Gamma
\subset \Gamma$ etc.

\begin{df}\label{dfgam}
Let $X$ be a normal noetherian scheme
and $U\subset X$ be a dense open subscheme.
Let $t\in U$,
$A\subsetneqq k(t)$ be a valuation ring 
of the residue field at $t$
and $T={\rm Spec}\, A\to X$
be a $U$-external morphism. 
Let
$\gamma \in (0,\infty)_{\Gamma_{\mathbf Q}}$
for $\Gamma_{\mathbf Q}=
\Gamma\otimes {\mathbf Q}$.
Let $Y$ be a quasi-finite flat scheme over $X$
such that 
$V=Y\times_XU\to U$ is \'etale.

We define a commutative
diagram 
\begin{equation}
\xymatrix{
F^{\infty}_T(Y/X)
\ar[d]_{\varphi_T^\gamma}
\ar[r]^{\varphi_T^{\gamma+}}
&
F^{\gamma+}_T(Y/X)
\ar[d]\ar[dl]
\\
F^\gamma_T(Y/X)
\ar[r]&
F^{0+}_T(Y/X)
}
\label{eqFg}
\end{equation}
as the inverse limit of
\begin{equation}
\xymatrix{
\bar Y'_{\bar x'}
\ar[d]_{\varphi_{\bar x'}^{R'}}
\ar[r]^{\hspace{-5mm}\varphi_{\bar x'}^{R'+}}
&
F^{R'+}_{\bar x'}(Y'/X')
\ar[d]\ar[dl]
\\
F^{R'}_{\bar x'}(Y'/X')
\ar[r]&
Y_{\bar x}}
\label{eqFR'}
\end{equation}
for commutative diagrams
{\rm (\ref{eqXT'})} satisfying
the conditions {\rm (i)-(iii)}.

We say that the ramification of
$Y$ over $X$ at $T$
is bounded by $\gamma$
(resp.\ by $\gamma+$) if
$F^{\infty}_T(Y/X)
\to F^\gamma_T(Y/X)$
(resp. 
$F^{\infty}_T(Y/X)
\to F^{\gamma+}_T(Y/X)$)
is an injection.
\end{df}

By Lemma \ref{lmcof}, 
the limit is a filtered limit.

\begin{lm}\label{lmRFg}
{\rm 1.}
There exist
a commutative diagram
{\rm (\ref{eqXT'})} satisfying
the conditions {\rm (i)-(iii)},
an effective Cartier divisor
$R'\subset X'$ satisfying $R'\cap U'=\varnothing$
and $x'\in R'$
such that
$Y'$ over $X'$ satisfies
{\rm (RF)} for $R'$
at the image $x'\in R'$ of the closed point of $T'$.

{\rm 2.}
For $x'\in R'\subset X'$
satisfying the condition in {\rm 1},
the canonical morphism from
{\rm (\ref{eqFg})} to
{\rm (\ref{eqFR'})} is an isomorphism.
The diagram
{\rm (\ref{eqFg})} is a diagram of finite sets.
\end{lm}

\proof{
1.
By Lemma \ref{lmcof}.1,
after replacing $X$ by
a normal scheme of finite type over $X$
if necessary,
we may assume that
there exist an effective
Cartier divisor $R\subset X$ such 
that $v(R)=\gamma$
and a closed immersion
$Y\to Q$ over $X$
to a smooth scheme $Q$ over $X$.
Applying Theorem \ref{thmRFT} to $Q^{(R)}\to X$
and taking the normalizations,
we obtain a morphism $X'\to X$
of finite type of normal noetherian 
schemes satisfying the following properties:
The morphism $X'\to X$
is the composition
of a blow-up $X^*\to X$ with center
supported in $X\sm U$
and a faithfully flat morphism
$X'\to X^*$ of finite type
such that
$U'=X'\times_XU\to U$ is \'etale.
The morphism $Q^{\prime (R')}\to X'$
is flat and reduced.
Hence $Y'$ over $X'$
satisfies the condition (RF) for $R'$.
The morphism $T\to X$ is lifted to
$T'\to X'$ by Lemma \ref{lmvalf}.

2.
By 1 and Lemma \ref{lmcof},
among commutative diagrams
(\ref{eqXT'})
those such that the base change 
$Y'=Y\times_XX'$ over $X'$
satisfies the condition (RF)
for $R'$ at $x'$
are cofinal.
Hence the assertion follows
from Lemma \ref{lmFunX}.
\qed }

\medskip

We study functoriality of
the construction of
$F^{\gamma}_T(Y/X)$
and
$F^{\gamma+}_T(Y/X)$.
We consider a commutative diagram
\begin{equation}
\begin{CD}
Y'@>>>X'@<<< T'@<<<\bar x'@<<< \bar t'\\
@VVV@VVV@VVV@VVV@VVV\\
Y@>>>X@<<< T@<<<\bar x@<<< \bar t,
\end{CD}
\label{eq5}
\end{equation}
dense open subschemes
$U\subset X$ and
$U'\subset U\times_XX'\subset X'$
and 
$\gamma\in (0,\infty)_{\Gamma_{\mathbf Q}}$ and
$\gamma'\in (0,\infty)_{\Gamma'_{\mathbf Q}}$
satisfying the following properties:

(i) $X'\to X$ is a morphism of normal
noetherian schemes.

(ii) $T={\rm Spec}\, A\to X$
and
$T'={\rm Spec}\, A'\to X'$
are $U$-external and
$U'$-external morphisms
for valuation rings
$A\subsetneqq K=k(t)$
and
$A'\subsetneqq K'=k(t')$
of the residue fields
at $t\in U$ and $t'\in U'$.
The morphism $T'\to T$ is faithfully flat.

(iii) $Y\to X$ and $Y'\to X'$
are quasi-finite morphisms
such that
$Y\times_XU\to U$ and
$Y'\times_{X'}U'\to U'$ are
\'etale.

(iv) $\gamma\leqq \gamma'$.

(v) $\bar x'\gets \bar t'$ is
a lifting of $\bar x\gets \bar t$.

\begin{lm}\label{lmfgm}
We keep the notation above.

{\rm 1.}
We have a commutative diagram
\begin{equation}
\begin{CD}
F^{\infty}_{T'}(Y'/X')
@>>>
F^{\gamma'+}_{T'}(Y'/X')
@>>>
F^{\gamma'}_{T'}(Y'/X')
@>>>
F^{0+}_{T'}(Y'/X')
\\
@VVV@VVV@VVV@VVV\\
F^{\infty}_T(Y/X)
@>>>
F^{\gamma+}_T(Y/X)
@>>>
F^{\gamma}_T(Y/X)
@>>>
F^{0+}_T(Y/X)
\end{CD}
\label{eqFung}
\end{equation}
of finite sets.
Further if $\gamma<\gamma'$,
we have an arrow
$$
F^{\gamma'}_T(Y'/X')
\to 
F^{\gamma+}_T(Y/X)$$
making the two triangles
obtained by dividing the middle
square commutative.

{\rm 2.}
If the left square in 
{\rm (\ref{eq5})} is cartesian
and if $\gamma=\gamma'$,
the vertical arrows
in {\rm (\ref{eqFung})}
are bijections.
\end{lm}

\proof{
By Lemma \ref{lmRFg}.1,
we may assume that
there exists an effective Cartier
divisor $R\subset X$ such that
$R\cap U=\varnothing$
and $v(R)=\gamma$ and that
$Y$ over $X$ satisfies the
condition (RF) for $R$.
Further by Lemma \ref{lmRFg}.1
and Lemma \ref{lmcof}.3,
we may assume that
there exists an effective Cartier
divisor $R'\subset X'$ such that
$R'\cap U'=\varnothing$,
$v'(R')=\gamma$ and
$R'\geqq R\times_XX'$ and that
$Y'$ over $X'$ satisfies the
condition (RF) for $R'$.
Then, by Lemma \ref{lmRFg}.2,
we may identify
$F^{\gamma}_T(Y/X)
=F^R_{\bar x}(Y/X)$,
$F^{\gamma+}_T(Y/X)
=F^{R+}_{\bar x}(Y/X)$
and
$F^{\gamma'}_{T'}(Y'/X')
=F^{R'}_{\bar x'}(Y'/X')$,
$F^{\gamma'+}_{T'}(Y'/X')
=F^{R'+}_{\bar x'}(Y'/X')$.

1.
The assertion now follows from
the functoriality of dilatation (\ref{eqFunDD'}).

2.
In the notation above,
we may further assume that
$R'=R\times_XX'$.
Hence the assertion follows from
Lemma \ref{lmFunX}.
\qed
}

\medskip

Let $T^h$ be the henselization
at the closed point
$x\in T$
and let $t^h\in T^h$ 
denote the generic point.
Then, the absolute Galois group
$D_T={\rm Gal}(\bar t/t^h)$
acts on the specialization
$\bar x\gets \bar t$
of geometric points of $T$.
Hence the commutative diagram
(\ref{eqFg}) admits a canonical
action of $D_T$.

\begin{thm}\label{thmval}
Let the notation be as in
Definition {\rm \ref{dfgam}}.
Then, there exist
an element $\beta_0
\in (0,\infty)_{
\Gamma_{\mathbf Q}}$ and
finite pairs 
$(\alpha_i,\beta_i)_{i\in I}$ of
elements of $[0,\beta_0]_{
\Gamma_{\mathbf Q}}$
satisfying the following properties
{\rm (i)-(iii)}:

{\rm (i)}
$[0,\beta_0]_{
\Gamma_{\mathbf Q}}
=\bigcup_{i\in I}
[\alpha_i,\beta_i]_{
\Gamma_{\mathbf Q}}$.

{\rm (ii)}  
For $\gamma>\beta_0$
(resp.\
$\gamma\geqq \beta_0$),
$F_T^{\gamma}(Y/X)
\gets F^{\infty}_T(Y/X)$
(resp.\ $F_T^{\gamma+}(Y/X)
\gets F^{\infty}_T(Y/X)$)
is an injection.

{\rm (iii)}
Let $i\in I$
and $\gamma\in
(\alpha_i,\beta_i)_{
\Gamma_{\mathbf Q}}$.
Then, 
$F^{\gamma}_T(Y/X)
\gets
F^{\beta_i}_T(Y/X)$
is an injection
and 
$F^{\alpha_i+}_T(Y/X)
\gets
F^{\gamma+}_T(Y/X)$
is an injection on the image of
$F^{\infty}_T(Y/X)$.
\end{thm}

\proof{
Since we may take base change,
we may assume that $Y\to X$ is finite
and that the normalization
$\bar Y\to X$ is finite \'etale.
Hence by Lemma \ref{lmb+},
we may assume that there
exists an effective Cartier divisor
$R\subset X$ such that
$R\cap U=\varnothing$ and
$\bar Y\to Y^{(R)}$
is a closed immersion.

Set $\beta_0=v(R)\in \Gamma$.
Then, by Lemma \ref{lmRFg},
after replacing $X$ if necessary,
we may assume that 
$Y$ over $X$ satisfies the condition
(RF) for $R$.
Since 
$\bar Y\to Y^{(R)}$
is a closed immersion
and $F^{\infty}_T(Y/X)
\to \bar Y_{\bar x}$ is
a bijection,
$\bar Y_{\bar x}=F^{\infty}_T(Y/X)\to
Y^{(R)}_{\bar x}=
F_T^{\beta_0+}(Y/X)$
is an injection.
For $\gamma>\beta_0$,
the composition
$F_T^{\beta_0+}(Y/X)
\gets
F_T^{\gamma}(Y/X)
\gets
F_T^{\gamma+}(Y/X)
\gets 
F^{\infty}_T(Y/X)$
is an injection.
Hence the condition (ii) is satisfied.

Let $Q$ be a smooth scheme over $X$
and let $Y\to Q$ be a closed immersion.
As in Example \ref{exCD},
we define a semi-stable curve
$C_R\to X$
by the effective Cartier divisor 
$R\subset X$.
Define an effective Cartier
divisor $\tilde R\subset C_R$
to be the exceptional divisor.
Applying Theorem \ref{thmsst}
to $(Q\times_XC_R)^{[\tilde R]}
\to C_R\to X$
and taking the normalizations,
we obtain a commutative
diagram
$$\begin{CD}
C_R@<<< C'\\
@VVV@VVV\\
X@<<< X'
\end{CD}$$
where 
$Y_{C'}=Y\times_XC'$ over $C'$
satisfies the condition (RF)
for $R'=\tilde R\times_{C_R}C'$
and 
$C'\to X'$ is a semi-stable
curve.

By Lemma \ref{lmvalf},
there exist a finite extension
$K'$ of $K$
and a valuation ring
$A'$ such that $A=A'\cap K$
and that 
$T\to X$ is lifted to
$T'={\rm Spec}\, A'\to X'$.
Let $x'\in X'$ denote
the image of the closed point of $T'$.
Further, for $\gamma\in [0,\beta_0]_
{\Gamma_{\mathbf Q}},$
after replacing $K'$ by a finite extension
if necessary, we may assume that
$\gamma$ is an element of
$[0,\beta_0]_{\Gamma'}.$


Let $I_1$ be the set of irreducible components
of the fiber $C'\times_{X'}x'$.
For $i\in I_1$,
let $C_i\subset C'\times_{X'}x'$
denote the corresponding
connected component.
Let $I_2$ denote the
set of singular points
of the fiber $C'\times_{X'}x'$.
For $i\in I_2$,
let $z_i\subset C'\times_{X'}x'$
denote the corresponding
singular point.
Set $I=I_1\amalg I_2$.

Since the assertion is \'etale local on
$X'$, we may assume that
for each $i\in I_1$,
there exists a section $s_i\colon X'\to C'$.
For $i\in I_1$,
set $\alpha_i=\beta_i
=v'(s_i^*R')\in \Gamma^{\prime+}$.
Since $\alpha_i=v'(\tilde R)$
for the composition $T'\to X'\to C'\to C_R$,
we have $\alpha_i\in [0,\beta_0]_{\Gamma_{\mathbf Q}}$.
For $i\in I_2$,
if $z_i$ is contained in
two irreducible components
$C_{i_1}$ and $C_{i_2}$
such that
$\alpha_{i_1}\leqq \alpha_{i_2}
\in \Gamma^{\prime+}$,
we define $\alpha_i
=\alpha_{i_1}
\leqq\beta_i=\alpha_{i_2}
\in [0,\beta_0]_{\Gamma_{\mathbf Q}}$.
If $z_i$ is contained in
a unique irreducible component
$C_{i_1}$,
we define $\alpha_i
=\beta_i=\alpha_{i_1}\in
[0,\beta_0]_{\Gamma_{\mathbf Q}}$.

We show that the condition (i) is satisfied.
Since $\alpha_i,\beta_i
\in [0,\beta_0]_{
\Gamma_{\mathbf Q}}$
for $i\in I$,
we have the inclusion
$[0,\beta_0]_{
\Gamma_{\mathbf Q}}
\supset
\bigcup_{i\in I}
[\alpha_i,\beta_i]_{
\Gamma_{\mathbf Q}}$.
Let $\gamma$ be an element
of $[0,\beta_0]_{
\Gamma_{\mathbf Q}}$.
Then, we may assume
$\gamma \in 
[0,\beta_0]_{\Gamma'}$.
Then, since $T'\to X$
has a lifting to $T'\to C_R$
such that $v'(\tilde R)=\gamma$
and since $C'\to C_R\times_XX'$
is proper and birational,
there exists a unique lifting
$T'\to C'$
of $T'\to C_R$
by the valuative criterion.
If the image of closed point
by $T'\to C'$
is contained in the smooth part
$C_i\cap C_x^{\prime{\rm sm}}$
of an irreducible component
$C_i\subset C'_x$ for $i\in I_1$,
then we have $\gamma=\alpha_i$.
If the image of closed point
by $T\to C'$
is the singular point $z_i\in C'_x$
for $i\in I_2$,
then we have $\gamma\in
[\alpha_i,\beta_i]_{
\Gamma_{\mathbf Q}}$
by Corollary \ref{corDiv}.
Thus, the condition (i)
is also satisfied.

We show that the condition (iii) is satisfied.
For $i\in I_1$
or $i\in I_2$ such that
$\alpha_i=\beta_i$,
there is nothing to prove.
Assume that $i\in I_2$
and 
$z_i$ is contained in
two irreducible components
$C_{i_1}$ and $C_{i_2}$
such that
$\alpha_i=\alpha_{i_1}
<\beta_i=\alpha_{i_2}
\in \Gamma^{\prime+}$
and let
$\gamma\in
(\alpha_i,\beta_i)_{
\Gamma_{\mathbf Q}}$.
Then, we may assume
$\gamma \in 
(\alpha_i,\beta_i)_{\Gamma'}$.
By Corollary \ref{corDiv},
after replacing $T'$ by an
extension if necessary,
we may take a morphism
$T'\to C'$ such that
the image of the closed point
$x'\in T'$ is $z_i$ and
$v'(R')=\gamma$.
Since $F^\gamma_T(Y/X)
=F^{R'}_{\bar z_i}(Y_{C'}/C')$
and $F^{\gamma+}_T(Y/X)
=F^{R'+}_{\bar z_i}(Y_{C'}/C')$
by Lemma \ref{lmfgm}.2
and Lemma \ref{lmRFg}.2,
the assertion follows from
Proposition \ref{prstD}.
\qed }

\medskip


We study some variants.

Let $X$ be a normal noetherian
scheme, $U$ be
a dense open subscheme
and let $V\to U$ be a finite
\'etale morphism.
We consider a cartesian diagram
\begin{equation}
\begin{CD}
Y'@<<<V\\
@VVV
\hspace{-10mm}
\square
\hspace{7mm}
@VVV\\
X'@<<<U\end{CD}
\label{eqXX'}
\end{equation}
of schemes of finite type over $X$
satisfying the following conditions:
The horizontal arrows are
dense open immersions,
$X'$ is normal,
$X'\to X$ is a proper birational
morphism inducing the identity
on $U$
and $Y'$ is finite flat over $X'$.

Let $A\subset K=k(t)$ be
a valuation ring of the residue field
at a point $t\in U$ and
$T={\rm Spec}\, A\to X$ be
a $U$-external morphism.
Let $x\in T$ denote the closed point
and $\bar x$ be a geometric point above $x$.
For 
$\gamma\in \Gamma_{{\mathbf Q},>0}$,
we define
\begin{equation}
\xymatrix{
F^{\infty}_T(V/U)
\ar[d]\ar[r]
&
F^{\gamma+}_T(V/U)
\ar[ld]\ar[d]
\\
F^{\gamma}_T(V/U)
\ar[r]&
F^{0+}_T(V/U)
}
\label{eqFVU}
\end{equation}
to be the inverse limit of
\begin{equation}
\xymatrix{
F^{\infty}_T(Y'/X')
\ar[d]\ar[r]
&
F^{\gamma+}_T(Y'/X')
\ar[ld]\ar[d]
\\
F^{\gamma}_T(Y'/X')
\ar[r]&
F^{0+}_T(Y'/X').
}
\label{eqFYX}
\end{equation}

Let $T_V$ denote the normalization
of $T$ in $V\times_XT$.
For $X'$ in (\ref{eqXX'}),
let $X'_T\subset X'$ denote
the reduced closed subscheme
supported on the closure of
$t\in U\subset X'$
and let $x'\in X'_T$
denote the image of the unique morphism
$T\to X'$ lifting $T\to X$.
Then,
since $A=\varprojlim_{X'\to X}
{\cal O}_{X'_T ,x}$,
we have 
$F^{0+}_T(V/U)
=T_V\times_T\bar x$.

\begin{lm}\label{lmTf}
Suppose that 
the normalization $T_V$ of
$T$ in $V\times_XT$ is
finite and flat over $T$.
Then, 
there exists a finite and flat
$Y'\to X'$
such that
$T_V=Y'\times_{X'}T$.
For such 
$Y'\to X'$,
the diagram
{\rm (\ref{eqFVU})}
is isomorphic to
{\rm (\ref{eqFYX})}.
\end{lm}

\proof{
Since $A=\varprojlim_{X'\to X}
{\cal O}_{X'_T ,x}$
in the notation above,
the existence of
finite flat $Y'\to X'$ 
such that
$T_V=Y'\times_{X'}T$ follows.
By the flattening theorem
\cite[Th\'eor\`eme (5.2.2)]{RG},
such $Y'\to X'$ are cofinal among
commutative diagrams (\ref{eqXX'}).
Hence the assertion follows
from Lemma \ref{lmfgm}.2.
\qed
}

\medskip
For a normal noetherian scheme $X$,
a formal ${\mathbf Q}$-linear combination
$R=\sum_ir_iD_i$
with positive coefficients $r_i\geqq 0$
of irreducible closed subsets 
$D_i$ of codimension
1 is called an effective 
${\mathbf Q}$-Cartier divisor
if a non-zero multiple is an effective
Cartier divisor.
The union $\bigcup_i D_i$
for $r_i>0$ is called the support of $R$.
For an open subset $U\subset X$,
if $U$ does not meet the support of $R$,
we write $R\cap U=\varnothing$
by abuse of notation.
For a $U$-external morphism
$T={\rm Spec}\, A\to X$,
the valuation $v(R)$ is defined
as an element of
$[0,\infty)_{\Gamma_{\mathbf Q}}$.

\begin{df}\label{dframR}
Let $X$ be a normal noetherian
scheme and
$U\subset X$ be a dense open subscheme.
Let $Y$ be a quasi-finite flat scheme over $X$
such that 
$V=Y\times_XU\to U$ is finite \'etale.
Let $R$ be an effective ${\mathbf Q}$-Cartier divisor of
$X$ such that $U\cap R$ is empty
and let $x\in X$ be a point 
contained in the support of $R$.

We say that the ramification of
$Y$ over $X$ is bounded by $R$
(resp.\ by $R+$) at $x$,
if for every $U$-external 
morphism $T\to X$,
the ramification of
$Y\to X$ is bounded by $v(R)$
(resp.\ by $v(R)+$)
in the sense of Definition 
{\rm \ref{dfgam}}.
\end{df}

\begin{lm}\label{lmR}
Let the notation be as in
Definition {\rm \ref{dframR}.}
Then, the following conditions 
{\rm (1), \rm (1$'$)}
and {\rm (2)} are equivalent:

{\rm (1)}
The ramification of
$Y\to X$ is bounded by $R$
(resp.\ by $R+$)
in the sense of Definition 
{\rm \ref{dframR}}.

{\rm (1$'$)}
The condition in Definition 
{\rm \ref{dframR}}
with $T$ restricted to be
a discrete valuation ring
is satisfied.

{\rm (2)}
For every morphism
$f\colon X'\to X$
of finite type of normal noetherian schemes 
such that 
$U'=U\times_XX'\to U$ is \'etale,
that $R'=f^*R$
is an effective Cartier divisor
and that 
$Y'=Y\times_XX'\to X'$ satisfies the condition 
{\rm (RF)} in Definition {\rm \ref{dfram}}
for $R'$,
the ramification of
$Y'\to X'$ is bounded by $R'$
(resp.\ by $R'+$)
at every point of $R'$
in the sense of Definition 
{\rm \ref{dfram}}.
\end{lm}

\proof{
(1$'$)$\Rightarrow$(2):
Let $X'\to X$ be as in (2)
and $x'\in R'$ be a point.
Let $X'_1\to X'$ be
the normalization of
the blow-up at the closure of
$x'$.
Then, the local ring
$A'={\cal O}_{X'_1,x'_1}$ 
at the generic point $x'_1$
of an irreducible component
of the inverse image of $x'$
is a discrete valuation ring.
The morphism
$T'={\rm Spec}\, A'\to X'_1\to X'$
is $U'$-external
and the image of the closed point
is $x'$.

For $\gamma'=v'(R')$,
by Lemma \ref{lmRFg}.2,
the commutative diagram (\ref{eqFR'})
is canonically identified with
$$
\xymatrix{
F^{\infty}_{T'}(Y'/X')
\ar[d]_{\varphi_{T'}^{\gamma'}}
\ar[r]^{\varphi_{T'}^{\gamma'+}}
&
F^{\gamma'+}_{T'}(Y'/X')
\ar[d]\ar[dl]
\\
F^{\gamma'}_{T'}(Y'/X')
\ar[r]&
Y'_{\bar x'}.
}
$$
Further, this commutative diagram
is canonically identified with
(\ref{eqFg})
for $\gamma=v(R)$
by Lemma \ref{lmfgm}.2.
Hence the assertion follows.

(2)$\Rightarrow$(1):
Let $T\to X$ be a $U$-external morphism
and $\gamma=v(R)$.
Then by Lemma \ref{lmRFg}.2,
the commutative diagram (\ref{eqFg})
is canonically identified with
 (\ref{eqFR'}).
Hence the assertion follows.

The implication
(1)$\Rightarrow$(1$'$)
is obvious.
\qed }

\begin{pr}\label{pretale}
Let the notation be as in
Definition {\rm \ref{dframR}}
and assume that the
ramification of $Y$ over
$X$ is bounded by $R+$.
Assume that $Y$ is
locally of complete intersection over $X$
and let
$$\begin{CD}
Q@<<< Y\\
@VVV
\hspace{-10mm}
\square
\hspace{7mm}
@VVV\\
P@<<<X
\end{CD}$$
be a cartesian diagram of
schemes over $X$
such that
$P$ and $Q$ are smooth over
$X$, the vertical arrows
are quasi-finite and flat
and the horizontal arrows
are closed immersions.

Let $X'$ be a normal noetherian scheme
over $X$ such that
$R'=R\times_XX'$ is
an effective Cartier divisor,
that
$Y'=Y\times_XX'$ over $X'$
satisfies the condition 
{\rm (RF)} for $R'$.

Then, the morphism $Q^{\prime (R')}
\to P^{\prime (R')}$
is \'etale on a neighborhood
of $Q^{\prime (R')}
\times_{X'}R'$.
\end{pr}

\proof{
First, we show that we may
assume that there exist
a closed subscheme
$Y'_0\subset Y'$
\'etale over $X'$,
an integer $n\geqq 1$
and an effective Cartier divisor
$D'_0\subset R'$ satisfying the following conditions:
We have an equality 
$R'_{Y'_0}=
R'_{Y'}$ of underlying sets.
Let ${\cal J}'_0
\subset {\cal O}_{R'_{Y'}}$
be the nilpotent ideal
defining
$R'_{Y'_0} \subset R'_{Y'}$.
Then, we have
${\cal J}_0^{\prime n}=0$ 
and $(n+1)D'_0= R'$.

Under the condition (RF),
the formation of $Q^{\prime (R')}
\to P^{\prime (R')}$
commutes with base change by
Lemma \ref{lmFunX} and
Example \ref{exQD}.1.
Since $Q^{\prime (R')}$
and $P^{\prime (R')}$ are flat over $X'$,
the \'etaleness of 
$Q^{\prime (R')}
\to P^{\prime (R')}$
is checked fiberwise.
Hence, we may take base change.
Let $x'\in R'$ be a point and
let $X''\to X'$ be the normalization
of the blow-up at the closure of $x'$.
Then, there exists a point $x''\in X''$
above $x'$ such that the local ring
${\cal O}_{X'',x''}$ is a discrete valuation ring.
Hence, by replacing $X'$ by
${\rm Spec}\, {\cal O}_{X'',x''}$,
we may assume that $X'$ is the spectrum
of a discrete valuation ring.

Then, we may assume that
$Y'\subset Q'$ is a union of sections
$X'\to Q'$.
There exists a disjoint union $Y'_0
\subset Y'$ of sections such that
we have an equality 
$R'_{Y'_0}=
R'_{Y'}$ of underlying sets.
Let $n\geqq 1$ be an integer
satisfying ${\cal J}_0^{\prime n}=0$ 
in the notation above.
After replacing $X'$ by a ramified covering if necessary,
there exists effective Cartier divisor $D'_0$ of $X'$
satisfying $(n+1)D'_0= R'$.

The finite morphism
$Y^{\prime (R')}\to X'$
is \'etale by Corollary \ref{coretaYD}.
Hence by the existence
of $Y'_0, D'_0$ and $n$
and by Lemma \ref{lmet},
the ${\cal O}_{Y^{\prime (R')}}$-module
${\cal O}_{Y^{\prime (R')}}
\otimes_{{\cal O}_Y} \Omega^1_{Y/X}$
is annihilated by ${\cal I}_{nD'_0}$.
Hence by Lemma \ref{lmetale},
there exists an
open neighborhood $W_1
\subset 
Q^{\prime (R')}$ of
$Q^{\prime (R')}
\times_{X'}R'$
such that 
$Q^{\prime (R')}\to
P^{\prime (R')}$
is \'etale on $W_1
\sm (Q^{\prime (R')}
\times_{X'}R')$.

The morphism
$Q^{\prime (R')}\to
P^{\prime (R')}$
is \'etale also on 
a neighborhood $W_2$ of
$Y^{\prime (R')}\subset
Q^{\prime (R')}$.
Since the vector bundle
$P^{\prime (R')}
\times_{X'}R'\to R'$ 
has irreducible fibers,
$W_2\subset Q^{\prime (R')}$
is dense in the fiber of
every point of $R'$ by
Proposition \ref{prCA}.1.
Hence the assertion follows from
 Lemma \ref{lmeta}.
\qed }

\subsection{Ramification groups}

\begin{thm}\label{thmgr}
Let $X$ be a connected
normal noetherian
scheme and $U\subset X$
be a dense open subscheme.
Let $G$ be a finite group,
$W\to U$ be a connected $G$-torsor 
and let $C$ be the category
of finite \'etale schemes over $U$
trivialized by $W$.
Assume that for every morphism
$V_1\to V_2$ of $C$,
the morphism $Y_1\to Y_2$
of normalizations of $X$
in $V_1$ and in $V_2$ is
locally of complete intersection.

Let $t\in U$ and
$T={\rm Spec}\, A\to X$
be a $U$-external morphism
for a valuation ring $A\subsetneqq 
K=k(t)$.
Let $\bar x$ (resp.\ $\bar t$)
be a geometric point
above the closed point $x$
(resp.\ the generic point $t$)
of $T$
and $\bar x\gets \bar t$
be a specialization.
Fix a lifting of $\bar x$
to the normalization $T_W$ of
$T$ in $W\times_XT$
and let $I_{\bar x}\subset G$
be the inertia group
at the image of the lifting of $\bar x$
to the normalization $Y_W$ of $X$ in $W$
by $T_W\to Y_W$.

For an object $V$ of $C$,
let $Y$ denote the normalization of
$X$ in $V$ and consider the fiber functor
sending $V$ to
$F^\infty_T(Y/X)$.

{\rm 1.}
There exist decreasing 
filtrations
$G^\gamma_T\supset
G^{\gamma+}_T$ of $G$ indexed by
$\gamma\in 
(0,\infty)_{\Gamma_{\mathbf Q}}$
such that, for every object $V$
of $C$, 
the canonical surjections
$F_T^{\infty}(Y/X)\to
F_T^{\gamma+}(Y/X)\to
F_T^{\gamma}(Y/X)$
induce bijections
\begin{equation}
G^{\gamma+}_T
\backslash F_T^{\infty}(Y/X)\to
F_T^{\gamma+}(Y/X),\qquad
G^\gamma_T
\backslash F_T^{\infty}(Y/X)\to
F_T^{\gamma}(Y/X).
\label{eqGg}
\end{equation}
For $I_{\bar x}=G_T^{0+}$,
the mapping
\begin{equation}
G^{0+}_T
\backslash F_T^{\infty}(Y/X)\to
F_T^{0+}(Y/X)
\label{eqG0}
\end{equation}
is a bijection.

{\rm 2.}
There exists a finite
increasing sequence 
$0=\alpha_0<\alpha_1<\cdots<\alpha_n$
of elements of 
$[0,\infty)_{\Gamma_{\mathbf Q}}$ such that we have
\begin{align}
G^{\alpha_{i-1}+}
=G^{\gamma}
=G^{\gamma+}
=G^{\alpha_i}
&\qquad
\text{for }
\gamma\in
(\alpha_{i-1},
\alpha_i)_{\Gamma_{\mathbf Q}},
1\leqq i\leqq n,
\\
G^{\alpha_n+}
=G^{\gamma}
=G^{\gamma+}
=1
&\qquad
\text{for }
\gamma\in
(\alpha_n,
\infty)_{\Gamma_{\mathbf Q}}.
\nonumber
\end{align}

{\rm 3.}
Let $D_T\subset G$ be the
decomposition group of $T$ in $W\times_XT$.
Then, 
$D_T$ normalizes 
$G^{\gamma}$ and
$G^{\gamma+}$.
\end{thm}

\proof{
1.
By Proposition \ref{prcoca}
and Lemma \ref{lmRFg}.2,
the diagram
\begin{equation}
\begin{CD}
F^\infty_T(Y'/X)
@>>>
F^{\gamma+}_T(Y'/X)
@>>>
F^\gamma_T(Y'/X)
@>>>
F^{0+}_T(Y'/X)
\\
@VVV
@VVV@VVV@VVV\\
F^\infty_T(Y/X)
@>>>
F^{\gamma+}_T(Y/X)
@>>>
F^\gamma_T(Y/X)
@>>>
F^{0+}_T(Y/X)
\end{CD}
\end{equation}
is a cocartesian diagram of
surjections.
Further the functors
$F^{\gamma}_T$ and
$F^{\gamma+}_T$ preserve
disjoint unions.
Hence by Proposition \ref{prGN},
we obtain filtrations
$(G^\gamma_T)_\gamma$
and
$(G^{\gamma+}_T)_\gamma$
indexed by
$\gamma
\in (0,\infty)_{\Gamma_{\mathbf Q}}$
characterized by the bijections (\ref{eqGg}).
For $\gamma=0$,
the bijection (\ref{eqG0})
follows from
$F_T^{0+}(Y/X)=Y_{\bar x}$.

2.
Since $C$ has only finitely
many connected objects
and 
$F^\infty_T(Y/X)
\to
F^{\gamma+}_T(Y/X)
\to
F^\gamma_T(Y/X)$ are surjections,
it follows from Theorem 
\ref{thmval}.

3.
Since the surjections
$F_T^{\infty}(Y/X)\to
F_T^{\gamma+}(Y/X)\to
F_T^{\gamma}(Y/X)$
are compatible with the actions
of $D_T\subset G$,
the subgroup 
$D_T\subset D_{\bar x}$ normalizes 
$G^{\gamma}$ and
$G^{\gamma+}$
by Corollary \ref{corGN}.
\qed }

\medskip

By the definition of the filtrations,
the ramification of
$Y/X$ at $T$ is bounded by
$\gamma$ (resp.\ by $\gamma+$)
if and only if the action of
$G_T^\gamma$ (resp.\ of $G_T^{\gamma+}$)
on 
$F_T^{\infty}(Y/X)$ is trivial.
By Corollary \ref{corqt},
the filtrations
$(G^{\gamma})$ and
$(G^{\gamma+})$
are compatible with quotients.
We have the following functoriality.
Let 
$$\begin{CD}
X'@<<< T'\\
@VVV@VVV\\
X@<<< T
\end{CD}$$
be a commutative diagram
of schemes.
Assume that
$X'\to X$ is a morphism
of normal connected noetherian
schemes and let
$U'\subset U\times_XX'
\subset X'$ be a dense open subscheme.
The horizontal arrows
$T\to X$ and $T'\to X'$
are $U$-external and
$U'$-external and
the vertical arrow $T'\to T$
is faithfully flat.
Let $W'$ be a connected 
$G'$-torsor over $U'$
for a finite group $G'$
and let $W'\to W$
be a morphism over $U'\to U$
compatible with a morphism
$G'\to G$ of finite groups.
Assume that $W'\to U'$ satisfies
the complete intersection property
as in Theorem \ref{thmgr}
and let 
$(G^{\prime\gamma'})$ and
$(G^{\prime\gamma'+})$
be the filtrations of $G'$
indexed by $\gamma'
\in (0,\infty)_{\Gamma'_{\mathbf Q}}.$
Then, for
$\gamma
\in (0,\infty)_{\Gamma_{\mathbf Q}},$
the morphism $G'\to G$
induces
\begin{equation}
G^{\prime\gamma}
\to G^{\gamma},
\quad
G^{\prime\gamma+}
\to G^{\gamma+}
\label{eqgrfn}
\end{equation}
by the functoriality Lemma \ref{lmfgm}.1.
\medskip

We consider a variant.
Let $A\subsetneqq K$
be a valuation ring
and $L$ be a finite Galois extension
of $K$ of Galois group $G$.
We define a filtration of $G$ by ramification 
groups under the following assumptions:
For every intermediate
extension $K\subset M\subset L$,
the normalization $A_M$ of
$A$ in $M$ is a valuation ring finite flat
and of complete intersection over $A$.
There exist an irreducible normal noetherian
scheme $X$ such that
$K$ is the residue field at the generic point $t$
and a morphism
$T={\rm Spec}\, A\to X$ extending
$t\to X$.

Let $T\to X$ be as above.
Then by Lemma \ref{lmcof}.1,
there exist a dense open subscheme 
$U\subset X$, 
a normal scheme
$X'$ of finite type over $X$
satisfying the following conditions:
The morphism
$U'=U\times_XX'\to U$
is an isomorphism.
The morphism
$T\to X$ is lifted to
$T\to X'$.
For every intermediate
extension $M$,
there exists a finite flat
scheme $Y'_M\to X'$
locally of complete intersection
such that
$U'\times_{X'}Y'_M\to U'$
is finite \'etale and 
$T\times_{X'}Y'_M={\rm Spec}\, A_M$.
Then applying Theorem \ref{thmgr},
we obtain filtrations
$(G_T^\gamma)$ and
$(G_T^{\gamma+})$ by normal subgroups
of $G=D_T$ indexed by
$(0,\infty)_{\Gamma_{\mathbf Q}}$.

\medskip


In the rest of the article,
we consider the case where
$X=T={\rm Spec}\,  {\cal O}_K$
for a complete discrete valuation ring
${\cal O}_K$.
For a finite Galois extension 
of the fraction field $K$
of Galois group $G$,
the decreasing filtrations 
$(G^r)_{r>0}$ and $(G^{r+})_{r\geqq 0}$ 
by normal subgroups
indexed by rational numbers 
are defined.

Let $L$ be a finite separable extension
of degree $n$ of $K$
and $Y={\rm Spec}\,  {\cal O}_L$
for the integer ring ${\cal O}_L$.
We recall the classical case where ${\cal O}_L$
is generated by one element
over ${\cal O}_K$,
using the Herbrand function.
Take a closed immersion
$Y={\rm Spec}\,  {\cal O}_L
\to Q={\mathbf A}^1_X=
{\rm Spec}\,  {\cal O}_K[T]$,
and let $P\in {\cal O}_K[T]$
be the monic polynomial
such that we have an isomorphism
${\cal O}_K[T]/(P)\to {\cal O}_L$.

Let $K'$ be a finite separable extension
containing the Galois closure of $L$
and $X'={\rm Spec}\, {\cal O}_{K'}$.
Let $v'\colon K'\to {\mathbf Q}
\cup\{\infty\}$ be the valuation
extending the normalized valuation of $K$.
Let $r>0$ be a rational number 
in the image of $v'$
and let $R'\subset X'$
be the effective Cartier divisor
such that $v'(R')=r$.
Let $Q'\supset Y'$ be the base change
of $Q\supset Y$ by $X'\to X$
and let $Q^{\prime (r)}=Q^{\prime (R')}$ 
denote the dilatation.
We compute $Q^{\prime (r)}$
using the Herbrand function
whose definition we 
briefly recall.

Decompose $P$ as
$P=\prod_{i=1}^n(T-a_i)$ in ${\cal O}_{K'}[T]$
and set $b_i=a_i-a_n\in {\cal O}_{K'}$.
Set $P(T_1+a_n)=\prod_{i=1}^n(T_1-b_i)
=T_1^n+c_1T_1^{n-1}+\cdots+c_{n-1}T_1$ 
in ${\cal O}_{K'}[T_1]$.
Changing the numbering if necessary,
we assume that the valuations 
$s_i=v'(b_i)\in {\mathbf Q}$ are increasing in $i$.
Note that the increasing sequence
$s_0=0\leqq 
s_1\leqq \cdots\leqq s_{n-1}<s_n=\infty$
is independent of the choice of $a_n$.
The valuation $v'(c_{n-1})=\sum_{k=1}^{n-1}s_k$
equals the valuation $v'(D_{L/K})$
of the different $D_{L/K}$.
It is further equal to
the length of the ${\cal O}_L$-module
$\Omega^1_{{\cal O}_L/{\cal O}_K}$
divided by the ramification index $e_{L/K}$
by \cite[Chap.\ III \S 7 Corollaire 2 \`a Proposition 11]{CL}.

The largest piecewise linear convex continuous
function $p\colon [0,n-1]
\to [0,v'(D_{L/K})]$ such 
that the graph is below
the points $(0,0)$ and $(k,{\rm ord}_Kc_k)$
for $k=1,\ldots,n-1$
is defined by
\begin{equation}
p(x)=\sum_{i=1}^{k-1}s_i+
s_k(x-k+1)
\label{eqNP}
\end{equation}
on $[k-1,k]$ for $k=1,\ldots,n-1$.
The graph of $p$ is
the {\em Newton polygon}
of the polynomial $P(T_1+a_n)$.
The {\em Herbrand function} 
$\varphi\colon [0,\infty)\to [0,\infty)$
is a piecewise linear concave continuous
function defined by
\begin{equation}
\varphi(s)=\sum_{i=1}^{n-1}
\min(s_i,s)+s.
\label{eqHer}
\end{equation}
We have 
\begin{equation}
\varphi(s)=\sum_{i=1}^{k-1}
(n-i+1)\cdot
(s_i-s_{i-1})+(n-k+1)\cdot (s-s_{k-1})
\label{eqphi}
\end{equation}
on $[s_{k-1},s_k)$ for $k=1,\ldots,n$.

\begin{ex}\label{exmono}
{\rm Let $s\in (s_{k-1},s_k]_{\mathbf Q}$,
$r=\varphi(s)$
and let $t$ be an element of 
a finite separable extension $K'$
of $K$
such that ${\rm ord}_Kt=s$.
By (\ref{eqphi})
and Example \ref{exQD},
$Q^{\prime(r)}$ is obtained
as an iterated dilatation
defined inductively by
$Q'_0=Q'$,
\begin{align*}
Q'_i&
=Q_{i-1}^{\prime ((n-(i-1))\cdot (s_i-s_{i-1}))}\quad
\text{ for }\ 0<i<k\ \text{ and}\\
Q^{\prime(r)}&
=
Q_{k-1}^{\prime ((n-(k-1))\cdot(s-s_{k-1}))}.
\end{align*}
Hence $Q^{\prime (r)}\to X'$ is smooth.
Let $C\subset Q^{\prime (r)}\times_{X'}x'$
be the connected component meeting 
the section $s'_n\colon X'\to Q^{\prime (r)}$
lifting $s_n\colon X\to Q$ defined by $T=a_n$.
Let $k$ be the smallest integer
$k=1,\ldots,n$
satisfying $s\leqq s_k$.
Then, 
${\rm Spec}\, {\cal O}_{K'}[T']$ for $T'=T_1/t$
is a neighborhood of $C\subset
Q^{\prime (r)}$.
Further on 
${\rm Spec}\, {\cal O}_{K'}[T']$,
the closed subscheme
$Y^{\prime (r)}
\subset Q^{\prime (r)}$
is defined by
$\prod_{i=k}^n(T'-b_i/t)$.

Consequently,
the surjection
$\bar Y_{\bar x}
=\{a_1,\ldots,a_n\}
\to F_X^r(Y/X)$
(resp.\ $\to F_X^{r+}(Y/X)$)
is given by the equivalence relation
$v'(a_i-a_j)\geqq s$
(resp.\ 
$v'(a_i-a_j)> s$).
In particular,
$r_{L/K}=\varphi(s_{n-1})=
v'(D_{L/K})+s_{n-1}$
is the unique rational number $r$
such that the ramification of
$Y$ over $X$ is bounded by $r+$
but not by $r$.
}
\end{ex}

We give a slightly simplified proof of
the proposition below
giving characterizations of
unramified extensions and
tamely ramified extensions.

\begin{lm}[{\rm \cite[Chap.\ III \S 7 Proposition 13]{CL},
\cite[Proposition A.3]{AS}}]\label{lmCL}
Let $L$ be a finite separable extension
of a complete discrete valuation field $K$.
Assume that ${\cal O}_L$
is generated by one element over ${\cal O}_K$
and let $r_{L/K}=\varphi(s_{n-1})=
v'(D_{L/K})+s_{n-1}$
be as in Example {\rm \ref{exmono}}.

{\rm 1.}
The following conditions are equivalent:

{\rm (1)}
$L$ is an unramified extension of $K$.

{\rm (2)}
$r_{L/K}=0$.

{\rm (3)}
$r_{L/K}<1$.

{\rm 2.}
The following conditions are equivalent:

{\rm (1)}
$L$ is a tamely ramified extension of $K$.

{\rm (2)}
$r_{L/K}=0$ or $1$.

{\rm (3)}
$r_{L/K}\leqq 1$.
\end{lm}

\proof{
By \cite[Proposition A.3]{AS},
we have
$v'(D_{L/K})\geqq 1-1/e_{L/K}$
and the equality holds if and only if
$e_{L/K}$ is tamely ramified.
We have $s_{n-1}\geqq 0$
and the equality holds if and only if
$L$ is unramified.
If $L$ is ramified,
we have $s_{n-1}\geqq 1/e_{L/K}$
and the equality holds
if and only if $L$ is tamely ramified.
The assertions follows from these.
\qed}

\begin{pr}[{\rm \cite[Proposition 6.8]{AS}}]\label{prtame}
Let $L$ be a finite separable extension of
a complete discrete valuation field $K$.

{\rm 1.}
The following conditions are equivalent:

{\rm (1)}
$L$ is an unramified extension of $K$.

{\rm (2)}
The ramification of $L$ over
$K$ is bounded by $1$.

{\rm 2.}
The following conditions are equivalent:

{\rm (1)}
$L$ is a tamely ramified extension of $K$.

{\rm (2)}
The ramification of $L$ over
$K$ is bounded by $1+$.
\end{pr}

\proof{
(1) $\Rightarrow$ (2):
Since ${\cal O}_L$ is generated by one element
over ${\cal O}_K$,
this follows from Example \ref{exmono}
and Lemma \ref{lmCL}.

(2) $\Rightarrow$ (1):
1.
Let $L$ be a finite separable extension
such that the ramification over
$K$ is bounded by $1$
and assume that $L$ was ramified over $K$.

Let $G$ be the Galois group
of a Galois closure of $L$ over $K$
and let $1\subsetneqq I
\subset G={\rm Gal}(L/K)$ 
be the wild inertia subgroup and the inertia subgroup.
By replacing $K$ and $L$ by
the subextensions corresponding to
$I$ and to a maximal subgroup
$H\subsetneqq I$, we may assume that
$L$ is a cyclic extension of prime degree 
since $I$ is solvable.

Then, either 
the ramification index $e_{L/K}$ is 1
and the residue extension is a
purely inseparable extension of degree $p$
or 
$L$ is totally ramified extension.
Hence ${\cal O}_L$ is generated by
one element and the assertion follows
from Example \ref{exmono}
and Lemma \ref{lmCL}.

2.
If the integer ring ${\cal O}_L$ is 
generated by one element over ${\cal O}_K$,
the assertion follows from Example \ref{exmono}
and Lemma \ref{lmCL}.
We prove the general case by reducing to
this case by contradiction.

Let $L$ be a finite separable extension
such that the ramification over
$K$ is bounded by $1+$
and assume that $L$ was wildly ramified over $K$.

Let $G$ be the Galois group
of a Galois closure of $L$ over $K$
and let $1\subsetneqq P\subset I
\subset G={\rm Gal}(L/K)$ 
be the wild inertia subgroup and the inertia subgroup.
By replacing $K$ and $L$ by
the subextensions corresponding to
$I$ and to a maximal subgroup
$H\subsetneqq P$, we may assume that
$[L:K]=mp$ for an integer $m$ prime to $p$.

Since an algebraic closure $\tilde F$ of the residue field $F$ of $K$
is a perfect closure
of the separable closure,
we may construct a henselian separable algebraic
extension $\tilde K$ of
ramification index 1 of
residue field $\tilde F$
as a limit
$\varinjlim K_\lambda$ of
finite separable
extensions of ramification index 1.
Since the composition 
$L\tilde K$ is a
totally ramified extension of
$\tilde K$,
there exists a finite
separable extension $K'=K_\lambda$
of ramification index 1
such that 
$L'=LK_\lambda$ is a
totally ramified extension of
$K'$.

By the functoriality 
(\ref{eqgrfn}),
the ramification of $L'$ over $K'$ is bounded by
$1+$.
Since $L'$ is totally ramified
over $K'$,
the integer ring ${\cal O}_{L'}$ is 
generated by one element over ${\cal O}_{K'}$.
Hence, 
$L'$ is tamely ramified
over $K'$ and we have $[L':K']=m$.

By construction, there exists a sequence
$K\subset K_0\subset K_1\subset
\cdots \subset K_n=K'$ 
such that $K_0$ is an
unramified extension of $K$
and that $K_i$ is an extension of $K_{i-1}$ of degree $p$
of ramification index $1$ 
with inseparable residue field extension
for each $i=1,\ldots, n$.
Since $[LK_0:K_0]=mp$,
we have $n>0$.
By taking the smallest such $n$,
we may assume
$[LK_{n-1}:K_{n-1}]=mp$. 

Further by the functoriality 
(\ref{eqgrfn}),
we may replace
$K$ and $L$ 
by $K_{n-1}$ and $LK_{n-1}$.
Hence, we may assume that
$[K':K]=p$ and $K'\subset L$.
Since $[K':K]=p$,
the integer ring ${\cal O}_{K'}$ is 
generated by one element over ${\cal O}_K$.
Since $K'\subset L$,
the ramification of $K'$ over
$K$ is bounded by $1+$.
Hence $K'$ is tamely ramified over $K$.
This contradicts to
that the residue field extension
of $K'$ over $K$ is inseparable.
\qed
}

School of Mathematical Sciences, University of Tokyo

Komaba, Meguro, Tokyo 153-8914 Japan

\end{document}